\documentclass[10pt,reqno]{amsart}
\usepackage{amssymb}
\usepackage{amsfonts}
\usepackage{textcomp}
\usepackage[centertags]{amsmath}
\usepackage{latexsym}
\usepackage{amsthm}
\usepackage{indentfirst}
\usepackage{fancyhdr}
\usepackage[english]{babel}
\usepackage[english]{babel}
\usepackage{color}
\usepackage{hyperref}
\usepackage[utf8]{inputenc}
\usepackage[T1]{fontenc}

\newcommand{\M}{\mathcal{M}}

\newtheorem{rmq}{Remark}[section]
\newenvironment{pre}[1][Proof]{\textbf{#1.} }{\ \rule{0.5em}{0.5em}}
\newtheorem{lem}{Lemma}[section]
\numberwithin{equation}{section}
\newtheorem{df}{Definition}[section]

\newtheorem{cor}{Corollary}[section]

\newtheorem{proposition}{Proposition}[section]
\newtheorem{thm}{Theorem}[section]
\usepackage{t1enc}
\title[]{\bf On eigenvalue problems for the $p(x)$-Laplacian}
\author[]{ Aboubacar Marcos \and  Janvier Soninhekpon} 
\date{}
\begin{document}
\maketitle

	\begin{abstract}
		This paper studies nonlinear eigenvalues problems  with a double non homogeneity  governed by the $p(x)$-Laplacian operator, under the  Dirichlet boundary condition on a bounded domain of $\mathbb{R}^N(N\geq2)$. According to the features of the nonlinearity (sublinear, superlinear) we use the Lagrange multiplier's method, the Ekeland's variational principle or the Mountain-Pass theorem to show that the spectrum includes  a continuous set of eigenvalues, which is in some contexts the whole set $\mathbb{R_+^{*}}$. Moreover, we show that the smallest  eigenvalue obtained  from the Lagrange multipliers is exactly  the first eigenvalue in the  Ljusternik-Schnirelman eigenvalues sequence and also provide sufficient conditions for multiplicity results.
	\end{abstract}	
	\section{Introduction}
	The search of eigenvalue  for non-homogeneous  problems  such as 
	\begin{equation}\label{0.0}
		D_1(\Omega)\left\lbrace\begin{array}{lll}
			{-\Delta_{p}}u=\lambda |u|^{q-2}u&\mbox{in $\Omega$ with }, q\neq p\\
			u=0 &\mbox{on $\partial\Omega$}
		\end{array}.\right.
	\end{equation}
can be reduced to study  the existence  of solution for problem
\begin{equation}\label{0.1}
		D'_1(\Omega)\left\lbrace\begin{array}{lll}
			{-\Delta_{p}}u= |u|^{q-2}u&\mbox{in $\Omega$}\\
			u=0 &\mbox{on $\partial\Omega$}
		\end{array}.\right.
	\end{equation}
	 since the parameter $\lambda$  can be scaled out by multiplying $u$ with a suitable real number. Of course when $u_1$ is a solution for problem  ($\ref{0.1}$), then  for any $\lambda \in (0, +\infty  ),   \lambda^{\frac{p-q}{q}}$ is an eigenvalue with eigenfunction $\lambda^{\frac{1}{p}}u_1.$ Accordingly, it is meaningless to seek  eigenvalues  for  problem  \eqref{0.0}  since then  the parameter $\lambda$ plays no role.
 When $p=2$, \eqref{0.1} becomes 
	\begin{equation}\label{0.11}
		D'_1(\Omega)\left\lbrace\begin{array}{lll}
			{-\Delta}u= |u|^{q-2}u&\mbox{in $\Omega$}\\
			u=0 &\mbox{on $\partial\Omega$}
		\end{array}.\right.
	\end{equation}
	and belongs to the celebrated family  equations of  Emden-Fowler which has been extensively 	studied in the literature.
	The question is quite different and complex  when dealing with eigenvalue problems involving a double non-homogeneity   such as 
	\begin{equation}\label{1.1}
		D_3(\Omega)\left\lbrace\begin{array}{lll}
			-\Delta_{p(x)}u=\lambda V(x)|u|^{q(x)-2}u&\mbox{in $\Omega$}\\
			u=0 &\mbox{on $\partial\Omega$}.
		\end{array}.\right.
	\end{equation}	
	The non-homogeneous operator $\Delta_{p(x)}$ is the so-called $p(x)$-Laplacian, $\Omega$ is a bounded domain in $\mathbb{R}^N (N\geq2)$ with smooth boundary,  $\lambda$ is a real parameter and $p(x) \neq q(x)$. We assume throughout the paper that $p$ and $q$ are continuous functions on $\overline{\Omega}$ and $V$ is a positive function in a generalized Lebesgue space $L^{s(x)}(\Omega)$.  
The $p(x)$- Laplacian operator appears in many contexts in physics, namely in  nonlinear electrorheological fluids and other phenomena related to image processing, elasticity and the flow in porous media; for a survey see (\cite{7,3,6,12,14,5}) and references therein.\\ 
	$\text{ }\text{ }\text{ }\text{ }\text{ }$
	When the variable exponents $p(.) = q(.)$, problem $(\ref{1.1})$ has been considered in several aspects and existence of eigenvalues and  some of their qualitative properties have been established (cf \cite{2,4,32, 15,3,33}) . 
	In the  particular case where  $p(x)=q(x)=\text{constant}$ and $V\equiv1$, An Lê proved in \cite{1} the existence of a non-decreasing sequence of nonnegative  Ljusternik-Schnirelman eigenvalues, and derived the properties of   simplicity and isolation of the principal eigenvalue.
For the case $p(x)=q(x)\neq\text{constant}$ and $V\equiv1$, Fan, Zhang and Zhao proved the existence of infinitely many eigenvalues sequence and  that $\sup\Lambda=+\infty$, where $\Lambda$ is the set of all nonnegative eigenvalues. They have also given some sufficient conditions under which $\inf\Lambda = 0$  or positive (see \cite{33}). In the same framework, an extension of the study to the whole space $\mathbb{R}^N$ has been carried out in \cite{53} by  N. Benouhiba.\\ 
$\text{ }\text{ }\text{ }\text{ }\text{ }$
When $q(x)\neq p(x)$, to the best of our knowledge, the pioneer work on  the eigenvalue problem \eqref{1.1} when $V = 1$, is  that  in \cite{41} where the authors prove the existence of a  continuous family of eigenvalues. They mainly suppose that 
\begin{equation}\label{assumption}
	1<\min_{x\in\overline{\Omega}}q(x)<\min_{x\in\overline{\Omega}}p(x)<\max_{x\in\overline{\Omega}}q(x)
\end{equation}  and show
	that there exists $\lambda^*>0$ such that any $\lambda\in(0,\lambda^*)$ is an eigenvalue for the problem $(\ref{1.1})$. One can clearly notice that  under their assumption the ranges of $p()$ and $q()$ can interfere as well.\\ 
$\text{ }\text{ }\text{ }\text{ }\text{ }$
 The present  paper considers the  problem $(\ref{1.1})$ in the case the weight $V$ is positive and in the both context of  sublinearity and superlinearity. The  purpose of the paper is to  investigate the structure of the spectrum of  problem  $(\ref{1.1})$   at the light of various methods of nonlinear analysis contributions. Roughly speaking, by means of  the Ekeland's variational principle , we show in a first part of the paper devoted to the sublinear case,  that  under  assumption (\ref{assumption}), problem $(\ref{1.1})$ admits   a continuous family of  eigenvalues in the neighborhood of $0$ for any  and particularly when the ranges of $p(.)$ and $q()$ do not interfere, we point out that the eigenvalues family is  exactly  the whole $\mathbb{R_+^{*}}$.  Moreover, we derive sufficient conditions for which  each eigenvalue of the continuous spectrum  admits, an infinitely countable (possibly unbounded) family of eigenfunctions.  Always in the context that the ranges of $p(.)$ and $q(.)$  do not interfere, we focused our investigation on  the eigenvalue problem  $(\ref{1.1})$ constrained to the sphere, using the Lagrange's multipliers method  and  next  on the Ljusternik-Schnirelman principle.  We derive from this part   that the smallest eigenvalue on the sphere provides by the  Lagrange's multipliers method  corresponds exactly to  the first eigenvalue of  the Ljusternik-Schnirelman  sequence. In the latter part, the superlinear (non coercive) problem has  been considered using  the Mountain- Pass theorem.  Globally, our work deals with sublinear and superlinear eigenvalue problems under many aspects of the exponent functions $p()$ and $q()$  and our approach involves new techniques contrasting with other treatments of  $(\ref{1.1})$ . \\
$\text{ }\text{ }\text{ }\text{ }\text{ }$
 This paper is organized as follows:
In section $\ref{section 2}$, we state some classical properties of the spaces $L^{p(x)}(\Omega)$ and $W^{1,p(x)}_0(\Omega)$. In section $\ref{main}$ we study the sublinear case  and we devote the section ${\ref{section7}}$ to the case where problem $(\ref{1.1})$ is superlinear.
	
	\section{General setting}\label{section 2}
	Let $C(\overline{\Omega})$ be the set of all continuous functions on $\overline{\Omega}$.
	Put $$C_+(\overline{\Omega})=\left\{h\in C(\overline{\Omega}) \text{ such that }h(x)>1 \text{ }\forall x \in \overline{\Omega}\right\}$$ For $p(.)\in C_+(\overline{\Omega})$, we put $$p^-=\min_{x\in\overline{\Omega}}p(x) \text{ }\text{ and }\text{ } p^+=\max_{x\in\overline{\Omega}}p(x)$$
	The variable exponent Lebesgue space $L^{p(x)}(\Omega)$ is defined by: $$L^{p(x)}(\Omega)=\left\{u: \Omega\longrightarrow\mathbb{R}\text{ measurable such that }\int_{\Omega}|u|^{p(x)}dx<\infty\right\}$$ with the norm $$\|u\|_{L^{p(x)}(\Omega)}=\|u\|_{p(x)}=\inf\left\{\lambda>0/\int_{\Omega}\left|\frac{u}{\lambda}\right|^{p(x)}dx\leq1\right\}$$\\
	\text{ }\text{ }\text{ }\text{ }\text{ }\text{ }\text{ }The variable exponent Sobolev space $W^{1,p(x)}(\Omega)$ is defined by: $$W^{1,p(x)}(\Omega)=\left\{u\in L^{p(x)}(\Omega)/\text{ } |\nabla u|\in L^{p(x)}(\Omega)\right\}$$ with the norm $$\|u\|_{W^{1,p(x)}(\Omega)}=\|u\|_{1,p(x)}=\inf\left\{\lambda>0/\int_{\Omega}\left(\left|\frac{u}{\lambda}\right|^{p(x)}+\left|\frac{\nabla u}{\lambda}\right|^{p(x)}\right)dx\leq1\right\}=\|u\|_{p(x)}+\|\nabla u\|_{p(x)}.$$ Its conjugate space is $L^{p'(x)}(\Omega)$ with $p'(x)=\frac{p(x)}{p(x)-1}$. 
	
	Define $W_0^{1,p(x)}(\Omega)$ as the closure of $C_0^{\infty}(\Omega)$ in $W^{1,p(x)}(\Omega)$. $W_0^{1,p(x)}(\Omega)$ is endowed with the norm $$\|u\|=\|\nabla u\|_{p(x)}$$  which is an equivalent norm on $W^{1,p(x)}(\Omega)$. For more information on the Lebesgue and Sobolev spaces with variables exponents, we refer to (\cite{7,18,8,9,6,10,5}).\\ Define the critical Sobolev exponent of $p$ by:
	$p^*(x)=\left\lbrace\begin{array}{lll}
		\frac{Np(x)}{N-p(x)} &\mbox{if $p(x)<N$}\\
		+\infty &\mbox{if $p(x)\geq N$}
	\end{array}\right.$ for all $x\in\Omega.$
	
	We have the following.
	\begin{proposition}\label{Prop2.4}\text{ }(cf \cite{11})
		\begin{itemize}
			\item[(1)] There exists $C_H>0$ such that for any $u\in L^{p(x)}(\Omega)$ and $v\in L^{p'(x)}(\Omega)$, we have the Holder's inequality: $$\left|\int_{\Omega}uv dx\right|\leq C_H\|u\|_{p(x)}\|v\|_{p'(x)}\leq 2\|u\|_{p(x)}\|v\|_{p'(x)}$$ where $C_H=\frac{1}{p^-}+\frac{1}{\left(p'\right)^-}$.
			\item[(2)] There is a constant $C>0$ such that for all $u\in W_0^{1,p(x)}(\Omega)$, $$\|u\|_{p(x)}\leq C\|\nabla u\|_{p(x)}:=C\|u\|$$
		\end{itemize}
	\end{proposition}
	\begin{proposition}\label{Prop2.1}
		(cf \cite{18,9,11,10}) If $p, q\in C(\overline{\Omega})$ and $1\leq q(x)<p^*(x)$ for all $x\in \overline{\Omega}$, then the embedding $W^{1,p(x)}(\Omega)\hookrightarrow L^{q(x)}(\Omega)$ is continuous and compact.
	\end{proposition}
	As consequences, we have:
	\begin{proposition}\label{Prop2.2}
		(cf \cite{18,9, 11,10}) If $p: \overline{\Omega}\longrightarrow (1,\infty)$ is continuous and $q: \Omega\longrightarrow (1,\infty)$ is a measurable function such that $p(x)\leq q(x)\leq p^*(x)$ a.e on $\overline{\Omega}$, then there is a continuous embedding $W^{1,p(x)}_0(\Omega)\hookrightarrow L^{q(x)}(\Omega)$.
	\end{proposition}
	\begin{proposition}\label{Prop2.5}
		If $p^+<\infty$, both spaces $\left(L^{p(x)},\|.\|_{p(x)}\right)$ and $\left(W^{1,p(x)}_0(\Omega), \|.\|\right)$ are separable, reflexive and uniformly convex Banach spaces.
	\end{proposition}
	Define the modular: $$\varphi_{p(x)}(u)=\int_{\Omega}|u|^{p(x)}dx.$$ Then we have the following properties:
	\begin{proposition}\label{Prop2.6}
		(cf \cite{13,11, 10}) For all $u, v\in L^{p(x)}(\Omega)$, we have :
		\begin{itemize}
			\item[(1)] $\|u\|_{p(x)}<1 (\text{resp. }=1, >1)\Leftrightarrow\varphi_{p(x)}(u)<1 (\text{resp. }=1, >1)$.
			\item[(2)] $\min\left(\|u\|_{p(x)}^{p^-}; \|u\|_{p(x)}^{p^+}\right)\leq\varphi_{p(x)}(u)\leq\max\left(\|u\|_{p(x)}^{p^-}; \|u\|_{p(x)}^{p^+}\right)$. Consequently, we have:\\$\left\lbrace\begin{array}{lll}
				\|u\|_{p(x)}^{p^-}\leq\varphi_{p(x)}(u)\leq\|u\|_{p(x)}^{p^+} &\mbox{if $\|u\|_{p(x)}>1$}\\
				\|u\|_{p(x)}^{p^+}\leq\varphi_{p(x)}(u)\leq\|u\|_{p(x)}^{p^-} &\mbox{if $\|u\|_{p(x)}\leq1$}
			\end{array}.\right.$
			\item[(3)] For $u_n, u\in L^{p(x)}(\Omega)$, we have: $u_n\rightarrow u\text{ if and only if }\varphi_{p(x)}(u_n-u)\rightarrow0$
			\item [(4)] For all $u, v\in L^{p(x)}(\Omega)$, $\varphi_{p(x)}(u+v)\leq 2^{p^+-1}\left(\varphi_{p(x)}(u)+\varphi_{p(x)}(v)\right)$.
		\end{itemize}
	\end{proposition}
	\begin{proposition}\label{Prop2.7}
		(cf \cite{13, 11,10})  Let $p$ and $q$ be measurable functions such that $p\in L^{\infty}(\Omega)$ and $1\leq p(x)q(x)\leq\infty$ for a.e $x\in\Omega$. Let $u\in L^{q(x)}(\Omega), u\neq0$. Then we have:$$\min\left(\|u\|_{p(x)q(x)}^{p^-}; \|u\|_{p(x)q(x)}^{p^+}\right)\leq\left||u|^{p(x)}\right|_{q(x)}\leq\max\left(\|u\|_{p(x)q(x)}^{p^-}; \|u\|_{p(x)q(x)}^{p^+}\right)$$ As consequence, we have: 
		$\left\lbrace\begin{array}{lll}
			\|u\|_{p(x)q(x)}^{p^-}\leq\left||u|^{p(x)}\right|_{q(x)}\leq\|u\|_{p(x)q(x)}^{p^+} &\mbox{if $\|u\|_{p(x)}>1$}\\
			\|u\|_{p(x)q(x)}^{p^+}\leq\left||u|^{p(x)}\right|_{q(x)}\leq\|u\|_{p(x)q(x)}^{p^-} &\mbox{if $\|u\|_{p(x)}\leq1$}
		\end{array}.\right.$
	\end{proposition}

	In the following, we put $X=W^{1,p(x)}_0(\Omega)$ with the norm $\|u\|=\|\nabla u\|_{p(x)}$. Define on $X$ the following functionals:
	\begin{equation}\label{2.1}
		F(u)=\int_{\Omega}\frac{V(x)}{q(x)}|u|^{q(x)}dx ,\text{ }\text{ }
		G(u)=\int_{\Omega}\frac{1}{p(x)}|\nabla u|^{p(x)}dx \text{, }\text{ }\phi(u)=\int_{\Omega}V(x)|u|^{q(x)}dx, \text{ }\text{ }\psi(u)=\int_{\Omega}|\nabla u|^{p(x)}dx
	\end{equation}

with $V \in L^{s(x)}(\Omega)$
	\begin{equation}
		I_\lambda(u)=\int_{\Omega}\frac{1}{p(x)}|\nabla u|^{p(x)}dx-\lambda\int_{\Omega}\frac{V(x)}{q(x)}|u|^{q(x)}dx
	\end{equation}
	It is well known that $F$ and $G$ belong to $C^1(X,\mathbb{R})$ (see \cite{20,24,19}) and one has for all $v\in X$:
	\begin{equation*}
		<F'(u),v>=\int_{\Omega}V(x)|u|^{q(x)-2}uvdx\text{ }\text{ }\text{ and }\text{ }\text{ }<G'(u),v>=\int_{\Omega}|\nabla u|^{p(x)-2}\nabla u\nabla vdx
	\end{equation*}

\noindent

 To go further in the setting of the functional framework we require the following 
\begin{equation}\label{A}
1<p(x), q(x) <N<s(x)  \textrm { and } s'(x)q(x)=\frac{s(x)q(x)}{s(x)-1}< p^*(x)   \textrm{ for any } x\in\overline{\Omega} 
\end{equation}
 From the assumptions, the embedding $W^{1,p(x)}(\Omega)\hookrightarrow L^{s'(x)q(x)}(\Omega)$ and $W^{1,p(x)}(\Omega)\hookrightarrow L^{q(x)}(\Omega)$ are continuous and compact.
	
	\begin{df}
		A pair $(\lambda,u)\in\mathbb{R}\times X$ is a weak solution of $(\ref{1.1})$ if:
		\begin{equation}\label{3.6}
			\int_{\Omega}|\nabla u|^{p(x)-2}\nabla u\nabla vdx=\lambda\int_{\Omega}V(x)|u|^{q(x)-2}uvdx \text{ }\text{ }\forall v\in X
		\end{equation}	
		Such a pair $(u,\lambda)\in X\times\mathbb{R}$ with $u$ non trivial is called an eigenpair, $\lambda$ is an eigenvalue and $u$ is called an associated eigenfunction.
	\end{df}
	It is well known that $u$ is a weak solution of problem $(\ref{1.1})$ on $X$ if and only if $u$ is a critical point of the energy functional $I_\lambda(u)=G(u)-\lambda F(u)$, that is:
	\begin{equation}\label{2.2}
		G'(u)-\lambda F'(u)=0
	\end{equation}
	To solve the eigenvalue problem $(\ref{2.2})$ the constrained variational method is usually employed (see \cite{43,44,2,45,46,1,47,38}). We take here $G$ as a constrained functional and $F$ as an objective functional. Let $\alpha>0$ and define: 
	\begin{equation}\label{Eq 2.7}
	\widetilde{M_\alpha} = \{u\in X; G(u)\leq \alpha>0\} 	
	\end{equation}
	and the boundary 
	\begin{equation}
	\partial\widetilde{M_\alpha} = M_\alpha=\{u\in X; G(u)=\alpha>0\}	
	\end{equation}
	It is well known that:
	\begin{itemize}
		\item [$\star$] $M_\alpha$ is a $C^1$-submanifold of $X$ with codimension one.
		\item [$\star$] $(u,\lambda)\in X\times\mathbb{R}$ satisfies relation $(\ref{3.6})$ if and only if $u$ is a critical point of $F$ with respect to $M_\alpha$ (see \cite{38}).
	\end{itemize}
	
	\begin{rmq}\label{rmq2.2} 
		For all $u\in M_\alpha$, we have $$\frac{\alpha p^-}{q^+F(u)}\leq\frac{\psi(u)}{\phi(u)}\leq\frac{p^+\alpha}{q^-F(u)}$$
	\end{rmq}

	The following propositions and definition  are useful  throughout this work.  Reader who is interested by any proofs  is invited to consult references therein.
	\begin{proposition}\label{Prop2.3}(cf \cite{53,4,26,17,3})
		The mapping $G$ is coercive, convex and sequentially weakly lower semi-continuous; that is $u_n\rightharpoonup u_0$ in $X$ implies $G(u_0)\leq\lim\limits_{n\rightarrow\infty}\inf G(u_n)$	
	\end{proposition}

	\begin{proposition}( cf \cite{3})\label{Prop 2.10}
		\begin{itemize}
			\item [(1)] $G':X\longrightarrow X^*$ is continuous, bounded and strictly monotone operator.
			\item [(2)] $G'$ is a mapping of type $(S_+)$; that is if $u_n\rightharpoonup u_0$ in $X$ and $\varlimsup_{n\rightarrow\infty}<G'(u_n)-G'(u_0);u_n-u_0>\leq0$ then $u_n\rightarrow u_0$ in $X$.
			\item [(3)] $G': X\longrightarrow X^*$ is a homeomorphism.
		\end{itemize}
	\end{proposition}

	\begin{proposition}( cf \cite{25})\label{Prop 2.9}
		Let $(X,d)$ be a complete metric space. Let $\Phi: X\longrightarrow\mathbb{R}\cup\{+\infty\}$ be lower semicontinuous and bounded below. Then given any $\varepsilon>0$ there exists $u_\varepsilon\in X$ such that:
		\begin{equation}
			\Phi(u_\varepsilon)\leq\inf_{u\in X}\Phi+\varepsilon
		\end{equation}
		and
		\begin{equation}
			\Phi(u_\varepsilon)<\Phi(u)+\epsilon d(u,u_\epsilon) \text{ }\forall u\in X\text{ with } u\neq u_\varepsilon
		\end{equation}
	\end{proposition}
	In the case the  functional $\Phi$ is of $C^1$ on Banach spaces, the Ekeland's variational principle takes the following version:
	\begin{cor}
		Let $J$ be a functional of class $C^1$ on a Banach space $X$, bounded below and $c=\inf_{u\in X}J(u)$. Then, given $\epsilon>0$, there is $u_\epsilon\in X$ such that:
		\begin{equation}\label{Eq 2.9}
			\left\lbrace\begin{array}{lll}
				c\leq J(u_\epsilon)\leq c+\epsilon\\
				\|DJ(u_\epsilon)\|_{X'}\leq\epsilon
			\end{array}.\right.
		\end{equation}
	\end{cor}

\begin{df}
 Let $X$ be a Banach space and $I:X \longrightarrow\mathbb{R}$ a continuously Fréchet differential functional. 
			
			We say that the functional $I$ satisfies the Palais-Smale condition if any sequence $(u_n)_n\subset X$ such that $I(u_n)$ is bounded and $I'(u_n)\rightarrow0$ in $X$ has convergent subsequence in $X$.
	\end{df}
	\begin{proposition}\label{Prop2.9}(cf \cite{36})
		Suppose that $X$ is a reflexive Banach space and let $K\subset X$ be a weakly closed subset of $X$. Suppose $J: X\longrightarrow\mathbb{R}\cup\{+\infty\}$ is coercive and sequentially-weakly lower semi-continuous on $K$ with respect to $X$. Then $J$ is bounded below and attains its minimum:
		\begin{equation*}
			\exists u\in K;\text{ } J(u)=\inf_{v\in K}J(v)=\min_{v\in K}J(v).
		\end{equation*}
	\end{proposition}
	Below we give the following proposition  and its proof which is adapted from the assumptions of our problem.
	\begin{proposition}\label{Prop 2.8}
		The functional $F$ is weakly-strongly continuous, that is, $u_n\rightharpoonup u$ in $X$ implies that $F(u_n)\longrightarrow F(u)$.
	\end{proposition}
	\begin{pre}
		Let $\left(u_n\right)_n\subset X$ be a sequence and $u\in X$ such that $u_n\rightharpoonup u$. We have:
		$$\left|F\left(u_n\right)-F(u)\right|\leq\frac{1}{q^-}\int_{\Omega}V(x)\left||u_n|^{q(x)}-|u|^{q(x)}\right|dx$$ By using the well known inequality
		\begin{equation*}
			\left||a|^p-|b|^p\right|\leq\gamma|a-b|\left(|a|+|b|\right)^{p-1}, \text{ }\text{ for all }p>1 \text{ and }(a,b)\in\mathbb{R}^2
		\end{equation*} where $\gamma$ is a positive constant, we obtain:
		\begin{equation}\label{Young1}
			\left|F\left(u_n\right)-F(u)\right| \leq \frac{\gamma}{q^-}\int_{\Omega}V(x)|u_n-u|\left(|u_n|+|u|\right)^{q(x)-1}dx
\end{equation}
and by applying two times  the Young inequality  to the right-hand side  term of the expression above, we get 
\begin{eqnarray*}
\int_{\Omega}V(x)|u_n-u|\left(|u_n|+|u|\right)^{q(x)-1}dx \leq&\\
 \frac{1}{s^-}\int_{\Omega}V^{s(x)}(x)dx + \frac{1}{q^- s'^-}\int_{\Omega}|u_n-u|^{s'(x)q(x)}dx + \frac{1}{q'^{-}s'^{-}} \int_{\Omega}\left(|u_n|+|u|\right)^{q(x)s'(x)}dx&
\end{eqnarray*}
Let's denote 
$a = \|V\|_{s(x)}, \quad b = \|u_n- u\|_{s'(x)q(x)}, \quad c = \||u_n| + |u|\|_{s'(x)q(x)}$\\
We have 
\begin{equation*}
\int_{\Omega}\frac{V(x)}{a}\frac{|u_n-u|}{b}\left(\frac{|u_n|+|u|}{c}\right)^{q(x)-1} dx \leq \frac{1}{s^-}+ \frac{1}{q^{-} s'^-} + \frac{1}{q'^{-}s'^{-}} .
\end{equation*}
Consequently, we have \\
$\int_{\Omega}V(x)|u_n-u|\left(|u_n|+|u|\right)^{q(x)-1} dx \leq$ $$ \frac{\gamma}{q^-} \left( \frac{1}{s^-}+ \frac{1}{q^{-} s'^-} + \frac{1}{q'^{-}s'^{-}}\right)  \|V\|_{s(x)}\|u_n- u\|_{s'(x)q(x)}\max( (\||u_n| + |u|\|)^{q^{+} -1}_{s'(x)q(x)}, (\||u_n| + |u|\|)^{q^{-} -1}_{s'(x)q(x)})$$
and then\\
$\left|F\left(u_n\right)-F(u)\right| \leq$ 
 $$C  \|V\|_{s(x)}\|u_n- u\|_{s'(x)q(x)}\max( (\||u_n| + |u|\|)^{q^{+} -1}_{s'(x)q(x)}, (\||u_n| + |u|\|)^{q^{-} -1}_{s'(x)q(x)}) $$
 where $C$ is a positive constant.
Since the embeddings $X\hookrightarrow L^{s'(x)q(x)}(\Omega)$  is compact, we have $\lim\limits_{n\rightarrow\infty}\|u_n-u\|_{s'(x)q(x)}=0$ and $\lim\limits_{n\rightarrow\infty}\left((\||u_n|\|+\||u|\|)_{s'(x)q(x)}\right)^{q^i-1}=2^{q^i-1}\|u\|_{s'(x)q(x)}^{q^i-1}$ for $i\in \{ -, +\}$. \\Hence $F(u_n)\rightarrow F(u)$ as $n\rightarrow\infty$.
	\end{pre}\text{ }\\

\section{Sublinear problem}\label{main}

	\subsection{Existence of a continuous family of eigenvalues}\label{section 3}
Let us consider for $V >0,$ the following Rayleigh quotients
	$${\mu}_*=\inf_{u\in X\smallsetminus\{0\}}\frac{\int_{\Omega}|\nabla u|^{p(x)}dx}{\int_{\Omega}V(x)|u|^{p(x)}dx}=
\inf_{u\in X\smallsetminus\{0\}}\frac{\psi(u)}{\phi(u)} ;\quad   \mu^*=\inf_{u\in X\smallsetminus\{0\}}\frac{\int_{\Omega}\frac{1}{p(x)}|\nabla u|^{p(x)}dx}{\int_{\Omega}\frac{1}{p(x)}V(x)|u|^{p(x)}dx}=\inf_{u\in X\smallsetminus\{0\}}\frac{G(u)}{F(u)} $$
\begin{equation*}
		\lambda_*=\inf_{u\in \widetilde{M_\alpha}}\frac{\int_{\Omega}|\nabla u|^{p(x)}dx}{\int_{\Omega}V(x)|u|^{q(x)}dx}=\inf_{u\in \widetilde{M_\alpha}}\frac{\psi(u)}{\phi(u)} \text{ }\text{ }\text{ and }\text{ }\lambda^*=\inf_{u\in \widetilde{M_\alpha}}\frac{\int_{\Omega}\frac{1}{p(x)}|\nabla u|^{p(x)}dx}{\int_{\Omega}\frac{V(x)}{q(x)}|u|^{q(x)}dx}=\inf_{u\in \widetilde{M_\alpha}}\frac{G(u)}{F(u)}
	\end{equation*}
We start this section by establishing the existence of a continuous family of eigenvalue in the case the ranges of $p(.)$ and $q(.)$ intercept. We recall that  assumptions $\eqref{A}$ are fulfilled, that is :
	\begin{equation*}\label{A1}
		1<p(x), q(x) <N<s(x)  \textrm { and } s'(x)q(x)=\frac{s(x)q(x)}{s(x)- 1}< p^*(x)   \textrm{ for any } x\in\overline{\Omega} .
	\end{equation*}
 For this aim, we consider the following eigenvalues sets:
	\begin{equation*}
		\Lambda=\{\lambda\in\mathbb{R}/\exists u\in X\smallsetminus\{0\} \text{ such that }(\lambda,u) \text{ is an eigenpair of } (\ref{1.1})\}
	\end{equation*}
	and for any $\alpha>0$,
	\begin{equation*}
		\widetilde{\Lambda_\alpha}=\{\lambda\in\mathbb{R}/\exists u\in \widetilde{M_\alpha} \text{ such that }(\lambda,u) \text{ is an eigenpair of } (\ref{1.1})\}.
	\end{equation*}
	Obviously, we have $\widetilde{\Lambda_\alpha}\subset\Lambda$. For $V>0$, we first point out the fact that if $(\lambda,u)$ is a solution of problem $(\ref{1.1})$, then $\lambda>0$. Indeed:
	\begin{eqnarray*}
		(\lambda, u) \text{ is a weak solution of }(\ref{1.1})&\iff & \int_{\Omega}|\nabla u|^{p(x)-2}\nabla u\nabla vdx=\lambda\int_{\Omega}V(x)|u|^{q(x)-2}uvdx \text{ }\text{ }\forall v\in X\\
		&\Rightarrow& \lambda=\frac{\int_{\Omega}|\nabla u|^{p(x)}dx}{\int_{\Omega}V(x)|u|^{q(x)}dx}\geq0 \text{ for }v=u.
	\end{eqnarray*}
	Suppose that $\lambda=0$. Thus $u$ is constant on $\overline{\Omega}$. This together with the fact that $u=0$ on $\partial\Omega$ gives $u\equiv0$ on $\overline{\Omega}$, which is a contradiction.\\
	Our goal in this first part of the paper is to investigate the existence of eigenvalues and corresponding eigenfunctions  in $\widetilde{M_\alpha}$ and on the whole $X.$
	\begin{thm}\label{Thm4.2'}
Assume that  assumptions $\eqref{A}$ are fulfilled
		\begin{itemize}
\item [(1)] If $\lambda$ is such that $0<\lambda<\lambda_*$ then $\lambda\notin\widetilde{\Lambda_\alpha}.$
\item [(2)] If  $q^- < p^- $ then for any $\alpha > 0$, there exists $\lambda_{\alpha} > 0 $  such that any $\lambda \in  (0, \lambda_{\alpha}) $  is an eigenvalue of problem $(\ref{1.1})$ with eigenfunction in $\widetilde{M_{\alpha}} .$ Moreover $\lambda_* = 0.$
		\end{itemize}
	\end{thm}
The proof will be carried out after the following  lemma. 
\begin{lem}\label{Lem 3.4}
		Let $\alpha>0$ ,  there exist $\lambda_\alpha>0$  such that for any $\lambda\in\left(0;\lambda_\alpha\right)$ we have:
\begin{equation}\label{eqlem 3.4}
		I_\lambda(u)\geq \frac{\alpha}{2}>0
\textrm {  for any } u\in M_\alpha.
\end{equation}
	\end{lem}

	\begin{pre}
		Let $\alpha>0$ and $u\in M_\alpha$. Trivially we have
		$ I_\lambda(u) = \alpha-\lambda\int_{\Omega}\frac{V(x)}{q(x)}|u|^{q(x)}dx.$ Thus we have 
		\begin{eqnarray*}
			I_\lambda(u) &=& \alpha-\lambda\int_{\Omega}\frac{V(x)}{q(x)}|u|^{q(x)}dx\\
			&\geq& \alpha-\frac{\lambda}{q^-}\int_{\Omega}V(x)|u|^{q(x)}dx\\
\end{eqnarray*}
and by Proposition $\ref{Prop2.6}$, we have:
		$$I_\lambda(u) \geq \alpha-\frac{\lambda}{q^-}C_HC^{q^\pm}\|V\|_{s(x)}\max(\|u\|^{q^-}, \|u\|^{q^+}).$$
On the other hand, since $u\in M_\alpha$, we have	
\begin{eqnarray}\label{Eq 3.3}
\|u\|&\leq &\max \left(\left(\alpha p^+\right)^{\frac{1}{p^-}}, \left(\alpha p^+\right)^{\frac{1}{p^+}}\right)	\text{ and next}\\
\max\left(\|u\|^{q^-}, \|u\|^{q^+}\right)&\leq &\max\left(\left(\alpha p^+\right)^{\frac{q^-}{p^-}}, \left(\alpha p^+\right)^{\frac{q^-}{p^+}}, \left(\alpha p^+\right)^{\frac{q^+}{p^-}}, \left(\alpha p^+\right)^{\frac{q^+}{p^+}}\right).
		\end{eqnarray}
	
		Consequently we derive
		\begin{equation*}
			 I_\lambda(u) 
			\geq \alpha-\frac{\lambda}{q^-}C_HC^{q^\pm}\|V\|_{s(x)}\max(\left(\alpha p^+\right)^{\frac{q^-}{p^-}}, \left(\alpha p^+\right)^{\frac{q^-}{p^+}}, \left(\alpha p^+\right)^{\frac{q^+}{p^-}}, \left(\alpha p^+\right)^{\frac{q^+}{p^+}}).
		\end{equation*}
		By the above inequality, we remark that if we define 
		\begin{equation}\label{Eq 3.6}
			\lambda_\alpha:=\frac{\alpha q^-}{2C_HC^{q^\pm}\|V\|_{s(x)}\max\left(\left(\alpha p^+\right)^{\frac{q^-}{p^-}}, \left(\alpha p^+\right)^{\frac{q^-}{p^+}}, \left(\alpha p^+\right)^{\frac{q^+}{p^-}}, \left(\alpha p^+\right)^{\frac{q^+}{p^+}}\right)},
		\end{equation}
		then for any $\lambda\in(0;\lambda_\alpha)$ and any $u\in M_\alpha$, we have 
		$$I_\lambda(u)\geq \frac{\alpha}{2}  >0.$$
	\end{pre}

\textbf{Proof of Theorem $\ref{Thm4.2'}$}
\begin{itemize}
			\item [(1)] Suppose that there exists $\lambda\in(0,\lambda_*)$ such that $\lambda\in \widetilde{\Lambda_\alpha}$. Thus there exists $u\in \widetilde{M_\alpha}$ such that $$\int_{\Omega}|\nabla u|^{p(x)-2}\nabla u\nabla vdx=\lambda\int_{\Omega}V(x)|u|^{q(x)-2}uvdx$$ for any $v\in X$. By taking $v=u$, we have $\int_{\Omega}|\nabla u|^{p(x)}dx=\lambda\int_{\Omega}V(x)|u|^{q(x)}dx$. Besides $$\lambda<\lambda_*=\inf_{u\in \widetilde{M_\alpha}}\frac{\int_{\Omega}|\nabla u|^{p(x)}dx}{\int_{\Omega}V(x)|u|^{q(x)}dx}\leq \frac{\int_{\Omega}|\nabla u|^{p(x)}dx}{\int_{\Omega}V(x)|u|^{q(x)}dx}.$$ 
Thus $$\int_{\Omega}|\nabla u|^{p(x)}dx>\lambda\int_{\Omega}V(x)|u|^{q(x)}dx=\int_{\Omega}|\nabla u|^{p(x)}dx.$$ This is a contradiction. 
\item[(2)]Let $\alpha>0$, $\lambda_\alpha$ be defined as in relation $(\ref{Eq 3.6})$ and $\lambda\in(0, \lambda_\alpha)$. Let consider the closed set $\widetilde{M_{\alpha}}$ as in $(\ref{Eq 2.7})$ and denote  $U_\alpha$   the interior set of the $\widetilde{M_{\alpha}}.$  $\widetilde{M_{\alpha}}$ endowed with the norm of $X$  is Banach subspace of $X$. Since $I_\lambda (0) =0 ,$   we deduce from Lemma $\ref{Lem 3.4}$ that $I_\lambda$ achieved its infimum in the interior set  of $\widetilde{M_{\alpha}}.$ 
Moreover
\begin{eqnarray}\label{bornitude1}
 I_\lambda (u) &\geq &-\frac{\lambda_{\alpha}}{q^-}C_HC^{q^\pm}\|V\|_{s(x)}\max(\|u\|^{q^-}, \|u\|^{q^+})\\
 &\geq &
- \frac{\lambda_{\alpha}}{q^-}C_HC^{q^\pm}\|V\|_{s(x)}\max\left(\left(\alpha p^+\right)^{\frac{q^-}{p^-}}, \left(\alpha p^+\right)^{\frac{q^-}{p^+}}, \left(\alpha p^+\right)^{\frac{q^+}{p^-}}, \left(\alpha p^+\right)^{\frac{q^+}{p^+}}\right)\nonumber
\end{eqnarray}

 that is $I_\lambda$ is bounded below  and then  
	$$-\infty < \inf_{\widetilde{M_\alpha}}I_\lambda=\inf_{U_\alpha}I_\lambda,$$ and consequently for any $\epsilon >0$, one can find $u_\epsilon \in U_\alpha $ such that
 $$I_\lambda(u_\epsilon)<\inf_{U_\alpha}I_\lambda +\epsilon = \inf_{\widetilde{M_\alpha}}I_\lambda+\epsilon.$$
And next, from the  Ekeland's  variational principle applied to the functional $I_\lambda:\widetilde{M}_{\alpha}\longrightarrow\mathbb{R}$, for any  $u \in \widetilde{M}_{\alpha}$  with  $u\neq u_\epsilon$, one has 
\begin{equation}\label{Ekeland1}
I_\lambda(u_\epsilon)<I_\lambda(u)+\epsilon\|u-u_\epsilon\| 
\end{equation}
Choose $u= u_\epsilon + t v$ for any $v \in \widetilde{M}_{\alpha}$   and $t>0$ small enough in \eqref{Ekeland1}, it then follows that
	$$\frac{I_\lambda(u_\epsilon+tv)-I_\lambda(u_\epsilon)}{t}\geq -\epsilon\|v\|;$$
	replacing $t$ by $-t$, $t$ chosen to be  negative in the inequality  above and letting $t\rightarrow0$, it follows  that $$\|I_\lambda'(u_\epsilon)\|_{X^*}\leq\epsilon.$$ 
We then deduce the existence of a sequence $(u_n)\subset U_\alpha$ such that
	\begin{equation}\label{Eq 3.21}
	-\infty< I_\lambda(u_n)\rightarrow \inf_{U_\alpha}I_\lambda = \underline{c} \text{ and }I_\lambda'(u_n)\rightarrow0 \text{ in }X^*	.
	\end{equation}
Obviously, $(u_n)_n$ is bounded and there exists $u_0\in X$ such that $u_n\rightharpoonup u_0$. It follows from sequential and  weak semicontinuity of $G$ ( cf Proposition $\ref{Prop2.3}$) that $G(u_0)\leq\alpha$; that is
 $u_0\in\widetilde{M}_{\alpha}$.
	On the other hand, we have:
	$$\langle G'(u_n),u_n-u_0\rangle=\langle I_\lambda'(u_n),u_n-u_0\rangle+\lambda\langle F'(u_n),u_n-u_0\rangle.$$
	Then using $(\ref{Eq 3.21})$ and the fact that $\langle F'(u_n),u_n-u_0\rangle\rightarrow0$ in Lemma $\ref{Lem 3.6}$, we deduce that 
		\begin{equation}\label{Eq 3.31}
		\langle G'(u_n),u_n-u_0\rangle\rightarrow 0,
		\end{equation}
	hence $(u_n)_n$ converges strongly to $u_0$ in X since the functional $G$ is of type $S_+$. 

 Since $I_\lambda\in C^1(X,\mathbb{R})$, we conclude that
	\begin{equation}\label{Eq 3.23}
	I_\lambda'(u_n)\rightarrow I_\lambda'(u_0) \text{ as } n\rightarrow\infty. 
	\end{equation}
	Relations $(\ref{Eq 3.21})$ and $(\ref{Eq 3.23})$ show that $$I_\lambda(u_0)=\underline{c} \text{ and }I'_\lambda(u_0)=0.$$
It remains to prove that $u_0\neq0.$ For this aim, it suffices to show that $c< 0.$ Indeed, since  $q^-<p^-$, we can choose  $\varepsilon >0$  such that $q^-+\varepsilon <p^-$. By the continuity of $q(.)$, we deduce the existence of an open set $\Omega_0\subset\Omega$ such that $q(x)\leq q^-+\varepsilon <p^-$ for all $x\in\Omega_0$.
Let $u_1\in M_\alpha$. It is obvious that $tu_1\in U_\alpha$ for any $t\in(0,1)$. So
		\begin{eqnarray*}
			I_\lambda(tu_1) &=& \int_{\Omega}\frac{t^{p(x)}}{p(x)}|\nabla u_1|^{p(x)}dx-\lambda\int_{\Omega}\frac{t^{q(x)}}{q(x)}V|v_0|^{q(x)}dx\\
			&\leq& t^{p^-}\alpha-\frac{\lambda}{q^+}\int_{\Omega_0}t^{q(x)}V(x)|u_1|^{q(x)}dx\\
			&\leq& t^{p^-}\alpha-\frac{\lambda t^{q^-+\varepsilon}}{q^+}\int_{\Omega_0}V(x)|u_1|^{q(x)}dx.
		\end{eqnarray*}
	Therefore 
	\begin{equation}\label{Eq 3.5}
		I_\lambda(tu_1)<0
	\end{equation}
	for $0<t<\delta^{\frac{1}{p^--q^--\varepsilon}}$ with
	$0<\delta<\min\left\{1;\frac{\lambda\int_{\Omega}V(x)|u_1|^{q(x)}dx}{\alpha q^+}\right\}.$\\
	Consequently, we get $$\inf_{U_\alpha}I_\lambda = c <0.$$
	So $u_0$ is a nontrivial weak solution for problem $(\ref{1.1})$ and thus any $\lambda\in(0;\lambda_\alpha)$ is an eigenvalue of problem $(\ref{1.1})$ with corresponding eigenfunction in $\widetilde{M}_{\alpha}$. From what precede, no eigenvalue lies in  $(0, \lambda_*)$ and hence  $\lambda_* = 0.$  $\blacksquare$
\end{itemize}

\begin{rmq}\label{remarq1}
	Expression of $\lambda_\alpha$ in \eqref{Eq 3.6} can be decoded to extract more information  on the size of  the eigenvalues set  with respect the exponent $p$ and $q$.\\  Indeed, suppose $\alpha p^+ \geq 1$  then 
\begin{equation}\label{Eq 3.6'}
			\lambda_\alpha:=\alpha^{1 -\frac{q^+}{p^-}}\frac{q^- \left( p^+\right)^{-\frac{q^+}{p^-}}}{2C_HC^{q^\pm}\|V\|_{s(x)}}
		\end{equation}
and hence $\lim\limits_{\alpha \rightarrow +\infty}\lambda_\alpha = +\infty \quad \text{when } \quad q^+ < p^-.$ On the other hand when $q^+ = p^-$, we observe that $$\lambda_\alpha = \lambda_{ p^-, q^+} := \frac{q^-}{2 p^+C_HC^{q^\pm}\|V\|_{s(x)}}$$  ceases to  depend on $\alpha .$ This fact will enable us to provide farther,  a multiplicity result on the eigenfunctions.  \\
In the case that  $\alpha p^+ < 1$ 
\begin{equation}\label{Eq 3.6''}
			\lambda_\alpha:=\alpha^{1 -\frac{q^-}{p^+}}\frac{q^- \left( p^+\right)^{-\frac{q^-}{p^+}}}{2C_HC^{q^\pm}\|V\|_{s(x)}}.
		\end{equation}
Hence, when $\alpha $ goes toward $0$,  $\lim\limits_{\alpha \rightarrow 0}\lambda_\alpha = 0$.\\
\end{rmq}

At the light of Remark $\ref{remarq1}$, we give in the following ,  the eigenvalues set $\widetilde{\Lambda_\alpha}$ when
problem $(\ref{1.1})$ is a sublinear problem. Roughly speaking we suppose that $1< q(x) < p(x) $ . 
\begin{cor}\label{Corol 1}
Assume that  assumptions  $\eqref{A}$ are fulfilled with  moreover  $q^+< p^-$,  then $\Lambda = (0,  +\infty)$, that is 
any $\lambda >0$ is an eigenvalue of problem $(\ref{1.1})$  on $X$ and hence  $\lambda_* =  \mu_* = 0$.
	\end{cor}
\begin{pre}	We know that $\Lambda \subset (0,  +\infty).$		
So let $\lambda \ge 0.$  Assuming that $q^+<p^-$, then from  Theorem $\ref{Thm4.2'}$, and for any $\alpha > 0,$  there is an eigenfunction $u_{\lambda} \in M_{\alpha}$ associated to $\lambda.$
$\lim\limits_{\alpha \rightarrow +\infty}\lambda_\alpha = +\infty$ and then there exists $\alpha$ big enough and  $\lambda_\alpha$ such that $\lambda \in (0, \lambda_\alpha)\subset \Lambda$. Thus $\Lambda = (0,  +\infty)$ with  eigenfunctions in $X.$\\
To conclude that  $\lambda_*=\mu_*=0$, we just have to notice that $\inf \widetilde{\Lambda_{\alpha}} = 0$ and $\lambda_*\geq \mu_*$ .
\end{pre}
\begin{cor}\label{Corol 2}
Assume that  assumptions  $\eqref{A}$ are fulfilled with $q^- < p^- $ and  $ q^+=  p^-$, then each  $\mu \in (0, \lambda_{q^+, p^-}), $ admits at least an infinitely countable  family of  eigenfunctions in $X$.
\end{cor}
\begin{pre}			
Let $\mu \in (0, \lambda_{ p^-, q^+}) $ and consider an increasing  sequence of  positive real numbers  $(\alpha_n )_{n>0}$  such that $\alpha_n p^+ \geq  1\quad  \forall  n. $ From Lemma $\ref{Lem 3.4}$ , there exists a positive real number   $\lambda_{\alpha_n}$  such that for any  $\lambda \in (0, \lambda_{\alpha_n}) $,  inequality $\eqref{eqlem 3.4}$ holds . But since $q^+ = p^-$ , $\lambda_{\alpha_ n} = \lambda_{ p^-, q^+}$ \quad for all $ n $,  then  $\mu \in (0, \lambda_{\alpha_n}) \quad \forall  n .$  Using inequality $\eqref{bornitude1}$ in the proof  of Theorem $\ref{Thm4.2'}$ along with the fact that  $\alpha_n p^+ \geq  1\quad  \forall  n $ and $q^+ = p^-$, we get 
\begin{eqnarray}\label{bornitude2}
 I_{\mu} (u) &\geq 
- \frac{\lambda_{ p^-, q^+}}{q^-}C_HC^{q^\pm}\|V\|_{s(x)}\max\left(\left(\alpha p^+\right)^{\frac{q^-}{p^-}}, \left(\alpha p^+\right)^{\frac{q^-}{p^+}}, \left(\alpha p^+\right)^{\frac{q^+}{p^-}}, \left(\alpha p^+\right)^{\frac{q^+}{p^+}}\right)\\
&\geq - \frac{\lambda_{ p^-, q^+}}{q^-}C_HC^{q^\pm}\|V\|_{s(x)}(\alpha_n p^+) \geq -\frac{\alpha_n}{2}.\nonumber
\end{eqnarray}
 So $I_{\mu}(u) $ is  bounded below  on $\widetilde{M}_{\alpha_n}$  and since $I_{\mu}(0) = 0 $  and  $I_{\mu}(u) \geq \alpha_n $ on $\widetilde{M}_{\alpha_n},$ it achieves its infimum in  $\widetilde{M}_{\alpha_n}$ interior .  Thus, processing closely to the idea developed in the the proof  of the  theorem, we  derive  on each $\widetilde{M}_{\alpha_n}$ some eigenfunctions associated to $\mu$ and consequently  when $n$ tends  toward $+\infty$, we can extract  a sequence of  eigenfunctions  $(u_n )_{n>0}$ belonging  to $X$  and having  the same eigenvalue  $\mu.$ Thus the proof is complete.
\end{pre}

\subsection{Eigenvalues problem constrained to a sphere }\label{section 4}

 The use of  the Lagrange multipliers to solve a constrained  eigenvalues problem like $(\ref{1.1})$  on the sphere $M_{\alpha}$ is reduced to  find a real number $\mu \in \mathbb{R}$ and   $u \in M_{\alpha} $ such that 
\begin{equation}\label{6.1}
			F'(u)=\mu G'(u), \text{ } \text{ }\text{ }\mu \in \mathbb{R}.
						\end{equation}
Accordingly $\lambda  = \frac{1}{\mu}$ will be an eigenvalue for problem $(\ref{1.1})$ corresponding to the eigenfunction $u$. Thus, we will point out in what follows that   problem $(\ref{1.1})$  admits an eigenvalue by means of the Lagrange multipliers method.
First of all let's denote
\begin{equation*}
		\nu_*=\inf_{u\in M_\alpha}\frac{\int_{\Omega}|\nabla u|^{p(x)}dx}{\int_{\Omega}V(x)|u|^{q(x)}dx}:=\inf_{u\in M_\alpha}\frac{\psi(u)}{\phi(u)} \text{ }\text{ }\text{ and }\text{ }\nu^*=\inf_{u\in M_\alpha}\frac{\int_{\Omega}\frac{1}{p(x)}|\nabla u|^{p(x)}dx}{\int_{\Omega}\frac{V(x)}{q(x)}|u|^{q(x)}dx}:=\inf_{u\in M_\alpha}\frac{G(u)}{F(u)}.
	\end{equation*}
\begin{equation*}
		\Lambda_\alpha=\{\lambda\in\mathbb{R}/\exists u\in M_\alpha \text{ such that }(\lambda,u) \text{ is an eigenpair of } (\ref{1.1})\}
	\end{equation*}

\noindent
\subsubsection{\bf Eigenvalue and Lagrange multiplier}\label{section 5}
\noindent

 The first result of this section is expressed as  follows.
\begin{thm}\label{Thm4.2}
 Consider that  assumption $\eqref{A}$ is fulfilled.  Then
		\begin{itemize}
\item [(1)]    $ \frac{q^-}{p^+}\nu_*\leq\nu^*\leq\frac{q^+}{p^-}\nu_*.$ and $\nu^* = 0$ \mbox{ if only if } $\nu_* = 0$ 
\item [(2)]$ \nu_* \neq 0 $  and if  $\lambda$ is such that $0<\lambda<\nu_*$ then $\lambda\notin\Lambda_\alpha$.
\item [(3)] If moreover $q^+<p^-$, then $\nu^*\notin\Lambda_\alpha$ and  there exists some $\lambda>\nu_*$ such that  $\lambda\in\Lambda_\alpha$. Moreover $\nu_*\notin\Lambda_\alpha$
		\end{itemize}
	\end{thm}
	The following lemmas are relevant for the proof of the theorem .
	
	\begin{lem}\label{Lem3.2}
		Let $\alpha>0$ be given. For any $u\in X\smallsetminus\{0\}$, there exists a unique $t>0$ such that $tu\in M_\alpha$.
	\end{lem}
	
	\begin{pre}
		Let $\alpha>0$ and $u\in X\smallsetminus\{0\}$ be given. Consider the function 
		\begin{eqnarray*}
			h: (0;+\infty) &\longrightarrow&(0;+\infty)\\
			t &\longmapsto& h(t)=G(tu)=\int_{\Omega}\frac{t^{p(x)}}{p(x)}|\nabla u|^{p(x)}dx
		\end{eqnarray*}
		Obviously the function $h$ is continuous and for any $t_1, t_2>0$ such that $t_1<t_2$, we have $h(t_1)<h(t_2)$; that is $h$ is strictly increasing. For $t\in(0,1)$, we have $h(t)\leq\frac{t^{p^-}}{p^-}\int_{\Omega}|\nabla u|^{p(x)}dx$ and thus $h(t)\longrightarrow0$ as $t\longrightarrow0$. For $t\in(1,\infty)$, $h(t)\geq\frac{t^{p^-}}{p^+}\int_{\Omega}|\nabla u|^{p(x)}dx$. Hence $h(t)\longrightarrow+\infty$ as $t\longrightarrow +\infty$. It follows that $h((0,\infty))=(0,\infty)$ and then the function $h$ is bijective. We deduce that for any $\alpha>0$, there exists a unique $t>0$ such that $h(t)=G(tu)=\alpha$; that is $tu\in M_\alpha$.
	\end{pre}
	\text{ }
	\\
	
	\begin{lem}\label{Lem 3.6}
		Let $\left(u_n\right)_n\subset X$ such that $u_n\rightharpoonup u$. Then 
		$$\lim\limits_{n\rightarrow\infty}\int_{\Omega}V(x)|u_n|^{q(x)-2}u_n(u_n-u)dx=0$$
	\end{lem}

	\begin{pre}
Proceeding similarly as in Proposition $\ref{Prop 2.8}$, we have 
		
			$\left|\int_{\Omega}V(x)|u_n|^{q(x)-2}u_n(u_n-u)dx\right| \leq$
			$$ \left(\frac{1}{s^-}+ \frac{1}{q^{-} s'^-} + \frac{1}{q'^{-}s'^{-}}\right)\|V\|_{s(x)}\max( \|u_n\|_{s'(x)q(x)}^{q^+-1}, \|u_n\|_{s'(x)q(x)}^{q^--1})\|u_n-u\|_{s'(x)q(x)}.$$
 From $u_n\rightharpoonup u_0$ and the compact embeddings $X\hookrightarrow L^{s'(x)q(x)}(\Omega)$, we have that $(|u_n)_n$  is bounded in  $L^{s'(x)q(x)}(\Omega)$ and $\|u_n-u\|_{s'(x)q(x)}\rightarrow0$.  The proof is complete.
	\end{pre}
	\text{ }\\
		
		Next, we move on to the proof of Theorem $\ref{Thm4.2}$\\
		\textbf{Proof of Theorem $\ref{Thm4.2}$}
		\begin{itemize}
\item [(1)]  First of all we observe that  $\nu^* = 0$ \mbox{ if only if } $\nu_* = 0 $ results easily from the inequality $ \frac{q^-}{p^+}\nu_*\leq\nu^*\leq\frac{q^+}{p^-}\nu_*.$\\
 Next, from Remark \eqref{rmq2.2} we have for all $u\in M_\alpha$,  $$\frac{\alpha p^-}{q^+F(u)}\leq\frac{\psi(u)}{\phi(u)}\leq\frac{p^+\alpha}{q^-F(u)}$$
And recalling the definitions of $\nu^*$  and $\nu_*$,  we get $ \frac{q^-}{p^+}\nu_*\leq\nu^*\leq\frac{q^+}{p^-}\nu_*.$ 
\item [(2)] Suppose  that $\nu^*=\inf\limits_{u\in M_\alpha\smallsetminus\{0\}}\frac{G(u)}{F(u)}=0$.  Thus there exists a sequence $(u_n)_n$ in $M_\alpha$ such that 
$\lim\limits_{n\rightarrow+\infty}\frac{G(u_n)}{F(u_n)}=\lim\limits_{n\rightarrow+\infty}\frac{\alpha}{F(u_n)}=0$; that is 
			\begin{equation}\label{4.1}
				\lim\limits_{n\rightarrow+\infty}F(u_n)=+\infty
			\end{equation}
			But since  the sequence $(u_n)_n$ belongs to  $M_\alpha$ , it is bounded  in $X$ and then converges weakly towards a function $u$  and because of the strong continuity  of  $F$ (cf Proposition $\ref{Prop 2.8}$), $\lim\limits_{n\rightarrow+\infty}F(u_n)= F(u)$. A contradiction with \eqref{4.1} and  then $\nu^*  > 0.$\\

Suppose that there exists $\lambda\in(0,\nu_*)$ such that $\lambda\in \Lambda_\alpha$. Arguing in a similar way as in the first assertion of  Theorem $\ref{Thm4.2'}$, we reach a contradiction and  hence there is no eigenvalue of problem $(\ref{1.1})$ in $(0,\nu_*).$
			
\item [(3)] Since $q^+<p^-$ and $\nu^*\leq\frac{q^+}{p^-}\nu_*$, we  get $\nu^*<\nu_*$ and since no eigenvalue belongs to  $(0, \nu_*)$, $\nu^*\notin\Lambda_\alpha$  . \\
			
	Let's prove that there exists  an eigenvalue $\lambda > \nu^* $   by proving  that  a minimizer of $I_\lambda$ on $M_\alpha$ is also a critical point. Since  
					\begin{equation}\label{AM1}
						\inf_{M_\alpha}I_\lambda(u)=\alpha-\lambda \sup_{M_\alpha}F(u),
					\end{equation}
					any minimizing sequence $(u_n)_n$ of $I_\lambda$ is a maximizing sequence of $F$. On the other hand  $(u_n)_n$ being in  ${M_\alpha}$ is bounded in $X$ and then, there exists $u_0 \in X$ such that  $(u_n)_n$ converges weakly to $u_0$. Because $G$ is sequentially weakly continuous and $F$ is strongly continuous, we have
					\begin{equation}\label{ball1}
						G(u_0) \leq \liminf\limits_{n\rightarrow +\infty} G(u_n)\leq \alpha,
					\end{equation}
					and 
					\begin{equation}\label{ball2}
						F(u_0) = \lim\limits_{n\rightarrow +\infty} F(u_n)=\sup_{M_\alpha}F(u).
					\end{equation}
					Clearly, from \eqref{ball1},  the maximum in \eqref{ball2} is achieved in the closed, convex and bounded set 
\begin{equation}\label{ball3}
\widetilde{M}_\alpha=\{u\in X; G(u)\leq \alpha\}. 
\end{equation}
But  $F$ is a convex function and  its maximum value  occurs on the boundary of   $\widetilde{M}_\alpha $, that is on  ${M}_\alpha $. Accordingly, the limit  $u_0 \in  {M}_\alpha $  and then 
					\begin{equation}
						I_\lambda(u_0)=\inf_{u\in M_\alpha}I_\lambda(u).
					\end{equation}
					
							Now,  let's show that $u_0$ is a critical point of $I_\lambda$  on  $M_\alpha$. Since $M_\alpha$ is not a vector space,  we  will consider some  small variations around  $u_0$ that lies on $M_\alpha$ so as in (see \cite{20}). 
							So, let $u\in X=W_0^{1,p(x)}(\Omega)$ be fixed and $\varepsilon>0$ small enough such that for any $s\in(-\varepsilon, \varepsilon)$, the function $s\mapsto u_0+su$ is not identically zero. From Lemma $\ref{Lem3.2}$, there exists a function $t: (-\varepsilon, \varepsilon) \longrightarrow (0,+\infty)$ such that 
							\begin{equation}\label{4.4}
								G(t(s)(u_0+su))=\int_{\Omega}\frac{[t(s)]^{p(x)}}{p(x)}|\nabla(u_0+su)|^{p(x)}dx=\alpha.
							\end{equation}
							As $u_0\in M_\alpha$, we deduce from relation $(\ref{4.4})$ that $t(0)=1$. On the other hand, we have 
							\begin{equation}
								\lim\limits_{s\rightarrow0}G\left(t(s)(u_0+su)\right)=G(u_0)=\alpha.
							\end{equation}
							Consequently, there is $\varepsilon>0$ small enough such that for any $s\in(-\varepsilon,\varepsilon)$ we have $t(s)(u_0+su)\in M_\alpha$. It follows that the map $s\longmapsto t(s)(u_0+su)$ defines a curve on $M_\alpha$ which passes through $u_0$ when $s=0$; that is $t(s)(u_0+su)\in M_\alpha$ for any $u\in X$ and $s\in(-\epsilon;\epsilon)$.
							
							One can easily see  that the function $t: s\in(-\epsilon,\epsilon)\longmapsto t(s)$ is derivable. Indeed, let's consider the map:
							\begin{eqnarray*}
								\tilde{g}: (-\epsilon,\epsilon)\times(0;+\infty) &\longrightarrow& \mathbb{R}\\
								(s,t) &\longmapsto& \tilde{g}(s,t)=\int_{\Omega}\frac{t^{p(x)}}{p(x)}|\nabla u_0+s\nabla u|^{p(x)}dx-\alpha.
							\end{eqnarray*}
							Obviously, $\tilde{g}(0, 1)=0$, $\tilde{g}$ is differentiable and for any $(s,t)\in(-\epsilon,\epsilon)\times(0;+\infty)$, we have: 
							$$\frac{\partial\tilde{g}}{\partial t}(s,t) =\int_{\Omega}t^{p(x)-1}|\nabla u_0+s\nabla u|^{p(x)}dx.$$
							We then deduce that 
							$$\frac{\partial\tilde{g}}{\partial t}(0, 1)=\int_{\Omega}|\nabla u_0|^{p(x)}dx\geq p^-\alpha \neq 0$$ and by means of $\eqref{eqlem 3.4}$ the  implicit functions theorem, for $\epsilon$ small enough, there exists a derivable function $\varphi:(-\epsilon,\epsilon)\longrightarrow(0;+\infty)$ such that $\forall (s,t)\in(-\epsilon,\epsilon)\times(0,+\infty), \tilde{g}(s,t)=0\Leftrightarrow t=\varphi(s)$ such that  $1=\varphi(0)$. Writing  $t(s)=\varphi(s)$, yields  a function $t: s\in(-\epsilon,\epsilon)\longmapsto t(s)$  derivable with $t(0) = 1$ .\\
							Moreover for any $s\in(-\epsilon,\epsilon)$ and $u\in X$ we have:
							
							$$\frac{\partial\tilde{g}}{\partial s}(s,t)=\int_{\Omega}t^{p(x)}|\nabla u_0+s\nabla u|^{p(x)-2}(\nabla u_0+s\nabla u)\nabla udx$$  and accordingly
							
							\begin{equation}
								t'(s)=-\frac{\frac{\partial\tilde{g}}{\partial s}(s,t)}{\frac{\partial\tilde{g}}{\partial t}(s,t)}=-\frac{\int_{\Omega}[t(s)]^{p(x)}|\nabla u_0+s\nabla u|^{p(x)-2}(\nabla u_0+s\nabla u)\nabla u dx}{\int_{\Omega}[t(s)]^{p(x)-1}|\nabla u_0+s\nabla u|^{p(x)}dx}.
							\end{equation}
							Put 
							\begin{equation}
								\gamma(s)=I_\lambda(t(s)(u_0+su)),
							\end{equation}
							of course $\gamma $ is derivable on  $(-\epsilon, \epsilon)$ and since $u_0$ is a minimal  point for $I_\lambda$,\quad  $s=0$ is a critical point for $\gamma$.\\ So, for any $s\in(-\epsilon;\epsilon)$, we have $$\gamma'(s) =\langle I_\lambda'(t(s)(u_0+su)); t'(s)(u_0+su)+t(s)u \rangle, \forall u \in X,$$ and hence :
							\begin{equation*}
								0=\gamma'(0) \Leftrightarrow \langle I_\lambda'(t(0)u_0); t'(0)u_0+t(0)u \rangle= \langle I_\lambda'(u_0); t'(0)u_0+u\rangle=0, \forall u \in X.
							\end{equation*}
							Recalling the expression of $t'$  we get 
							$t'(0)=-\frac{\int_{\Omega}|\nabla u_0|^{p(x)-2}\nabla u_0\nabla udx}{\int_{\Omega}|\nabla u_0|^{p(x)}dx}=-\frac{\langle G'(u_0),u\rangle}{\int_{\Omega}|\nabla u_0|^{p(x)}dx}$ and then 
							\begin{equation*}
								\langle I_\lambda'(u_0),u \rangle=-t'(0)\langle I'_\lambda(u_0),u_0\rangle, \forall u \in X
							\end{equation*}
							that is
							\begin{equation*}
								\langle G'(u_0),u\rangle-\lambda\langle F'(u_0),u\rangle=\frac{\langle G'(u_0),u\rangle}{\int_{\Omega}|\nabla u_0|^{p(x)}dx}\left(\int_{\Omega}|\nabla u_0|^{p(x)}dx-\lambda\int_{\Omega}V(x)|u_0|^{q(x)}dx\right), \forall u \in X.
							\end{equation*}
							Hence
							\begin{equation*}
								\langle G'(u_0),u\rangle-\frac{\int_{\Omega}|\nabla u_0|^{p(x)}dx}{\int_{\Omega}V(x)|u_0|^{q(x)}dx}\langle F'(u_0),u\rangle=0,  \forall u \in X
							\end{equation*}
							and then
							\begin{equation*}
								\langle I'_\lambda(u_0),u\rangle=0, \quad \forall u \in X \text{ with } \lambda=\frac{\int_{\Omega}|\nabla u_0|^{p(x)}dx}{\int_{\Omega}V(x)|u_0|^{q(x)}dx}.
							\end{equation*}
							
							Hence $\lambda>\nu_*$  is an eigenvalue of problem $(\ref{1.1})$ with its corresponding eigenfunction $u_0$ in $M_\alpha$ and is the smallest one constrained to the sphere $M_\alpha$ .
							
							We next show that $\nu_*\notin \Lambda_\alpha$. 
							Suppose by contradiction that here exists $u_*\in M_\alpha$ such that 
							\begin{equation}\label{Eq3.23}
								\langle I'_{\nu_*}(u_*),v\rangle=0 \text{ for any } v\in X.
							\end{equation}
							By taking $v=u_*$ in equation $(\ref{Eq3.23})$, we obtain:
							\begin{equation}
								\nu_*=\frac{\psi(u_*)}{\phi(u_*)}
							\end{equation}
							Obviously, we have $(1+s)t(s)u_*\in M_\alpha$ for any $s\in(-\epsilon;\epsilon)$ and thus
							\begin{equation}\label{Eq 3.25}
								\frac{\psi(u_*)}{\phi(u_*)}=\nu_*\leq\frac{\psi((1+s)t(s)u_*)}{\phi((1+s)t(s)u_*)} \text{ for any }s\in(-\epsilon;\epsilon)
							\end{equation}
							We will show that there are  some $s_0\in(-\epsilon, \epsilon)$ such that $0<[(1+s_0)t(s_0)]^{p^--q^+}<1$. For this purpose, we define the function
							\begin{equation}\label{Eq3.17}
								g(s)=G(t(s)(u_*+su))=\alpha\text{ }\text{ } \forall (u,s)\in X\times(-\epsilon,\epsilon)
							\end{equation}
							The function $g$ is derivable and we have for all $u\in X$ and $s\in(-\epsilon, \epsilon)$:
							\begin{equation}\label{Eq3.7}
								0=g'(s)=\langle G'(t(s)(u_*+su));t'(s)(u_*+su)+t(s)u\rangle
							\end{equation}
							Let $\theta>1$. For $u=\theta u_*$ in relation $(\ref{Eq3.7})$, we have:
							\begin{eqnarray*}
								0=g'(s) &=& \langle G'((1+\theta s)t(s)u_*);\left((1+\theta s)t'(s)+\theta t(s)\right)u_*\rangle\\
								&=& \left((1+\theta s)t'(s)+\theta t(s)\right)\langle G'((1+\theta s)t(s)u_*);u_*\rangle\\
								&=&\left((1+\theta s)t'(s)+\theta t(s)\right)\int_{\Omega}((1+\theta s)t(s))^{p(x)-1}|\nabla u_*|^{p(x)}dx
							\end{eqnarray*}
							So  for any $s\in(0 , \epsilon)$ with $\epsilon>0$ small enough, we have:
							\begin{eqnarray*}
								\int_{\Omega}((1+\theta s)t(s))^{p(x)-1}|\nabla u_*|^{p(x)}dx &\geq& ((1+\theta s)t(s))^{p^i-1 }\int_{\Omega}|\nabla u_*|^{p(x)}dx\\
								&\geq& ((1+\theta s)t(s))^{p^i-1 }p^+\alpha\\
								&>&0
							\end{eqnarray*}
							We then deduce that, for any $s\in(0,\epsilon)$:
							\begin{eqnarray*}
								g'(s)=0 &\Leftrightarrow& (1+\theta s)t'(s)+\theta t(s)=0\\
								&\Leftrightarrow& t'(s)+\frac{\theta}{1+\theta s}t(s)=0\\
								&\Leftrightarrow& t(s)=\tau e^{-\ln(1+\theta s)} \text{ with }\tau\in\mathbb{R}\\
								&\Leftrightarrow& t(s)=\frac{1}{1+\theta s} \text{ since } t(0)=1\\
								&\Leftrightarrow& (1+s)t(s)=\frac{1+s}{1+\theta s}.
							\end{eqnarray*}
							For some $s_0\in(0,\epsilon)$ and any $\theta>1$, we have $0<(1+s_0)t(s_0)=\frac{1+s_0}{1+\theta s_0}<1$. Since $p^->q^+$, we deduce that $0<[(1+s_0)t(s_0)]^{p^--q^+}<1$ . We then deduce from relation $(\ref{Eq 3.25})$ that
	\begin{equation}\label{Eq 3.29}
	\frac{\psi(u_*)}{\phi(u_*)}\leq\frac{\psi((1+s_0)t(s_0)u_*)}{\phi((1+s_0)t(s_0)u_*)}\leq [(1+s_0)t(s_0)]^{p^--q^+}\frac{\psi(u_*)}{\phi(u_*)}<\frac{\psi(u_*)}{\phi(u_*)}.
	\end{equation}
	By combining relations $(\ref{Eq 3.25})$ and $(\ref{Eq 3.29})$ we have
	$\frac{\psi(u_*)}{\phi(u_*)}<\frac{\psi(u_*)}{\phi(u_*)}$ 
	which is a contradiction.
							
	In short, we prove that if $\nu_*$ is an eigenvalue of problem $(\ref{1.1})$ with corresponding eigenfunction $u_*\in M_\alpha$, then $\langle I'_{\nu_*}(u_*);\theta u_*\rangle\neq0$ for any $\theta>1$. Hence $\nu_*\notin\Lambda_\alpha$ if $p^->q^+$.							
\end{itemize}

\noindent
										
	\subsubsection{\bf Existence of a Ljusternik-Schnirelman eigenvalues sequence}\label{section 6}
\noindent

	The most popular characterization of eigenvalues for nonlinear operator is certainly  due to the Ljusternik-Schnirelman principle. Thus, many  results on the  Ljusternik-Schnirelman characterizations exist in the literature when $ p= q $ is a  constant or  when $p(.) = q(.)$  under various assumptions and  boundary conditions (\cite{2, 4,  45'}). Here we are interested in the case $p(.) \neq q(.)$ and particularly when assumptions $\eqref{A}$ are satisfied and  $ q^+  < p^-$  that is when the ranges of $p$ and $q$ do not interfere. By means of a version of the Ljusternik-Schnirelman (L-S) principle (see \cite{44}) we  derive  the existence of  a sequence of eigenvalues for problem $(\ref{1.1})$ in $\M_{\alpha}\subset W^{1,p(x)}(\Omega)$ and moreover,  we establish a relationship between the  smallest eigenvalue yields by the multiplier of Lagrange method  in Theorem $\ref{Thm4.2}$ and the  first eigenvalue  in the (L-S) sequence .\\

	Let $X$ be a real reflexive Banach space and $F$, $G$ some functionals in $X$ as above. 
						\\We assume that:
						\newline (H1): $F, G: X\longrightarrow\mathbb{R}$ are even functionals and that $F, G \in C^1(X,\mathbb{R})$ with $F(0)=G(0)=0$.
						\newline (H2): $F'$ is strongly continuous (i.e $u_n\rightharpoonup u$ in $X$ implies $F'(u_n)\longrightarrow F'(u)$) and $ \langle F'(u),u\rangle=0, u\in \overline{coM_{\alpha}}$ implies $F(u)=0$ where $\overline{coM_{\alpha}}$ is the closed convex hull of $M_{\alpha}$ where $M_{\alpha}$ is as in the previous sections.
						\newline (H3): $G'$ is continuous, bounded and satisfies condition $S_0$; i.e as $n\rightarrow\infty$, $u_n\rightharpoonup u$; $G'(u_n)\rightharpoonup v$; $\langle G'(u_n), u_n\rangle\rightarrow \langle v,u\rangle$ implies $u_n\rightarrow u$.
						\newline (H4): The level set $M_{\alpha}$ is bounded and $u\neq 0$ implies $\langle G'(u),u\rangle>0$, $\lim\limits_{t\rightarrow\infty}G(tu)=+\infty$ and $\inf\limits_{u\in M_{\alpha}}\langle G'(u),u\rangle>0$.
						
						Let
						\begin{equation*}
							\Sigma_{(n,\alpha)}=\{H\subset M_{\alpha}; H \text{ is compact}, -H=H \text{ and } \gamma(H)\geq n\}
						\end{equation*}	
						where $\gamma(H)$ denotes the genus of $H$, i.e $\gamma(H):=\inf\left\{k\in \mathbb{N}/\exists h: H\longrightarrow\mathbb{R}^k-\{0\} \text{ such that } h \text{ is continuous and odd }\right\}$.\\
						Let's define the following  (L-S) sequence 
						\begin{equation}
							a_{(n,\alpha)}=\left\lbrace\begin{array}{ll}
								\sup\limits_{H\in\Sigma_{(n,\alpha)}}\inf\limits_{u\in H}F(u)&\mbox{if $\Sigma_{(n,\alpha)}\neq\phi$}\\
								0 &\mbox{if $\Sigma_{(n,\alpha)}=\phi$}	
							\end{array}.\right.
						\end{equation}

						\begin{equation}
							\chi_{\alpha}=\left\lbrace\begin{array}{ll}
								\sup\{n\in\mathbb{N}; a_{(n,\alpha)}>0\}&\mbox{if $a_{(1,\alpha)}>0$}\\
								0 &\mbox{if $a_{(1,\alpha)}=0$}
							\end{array}.\right.
						\end{equation}
						\\
						We suppose that $q^+ < p^-$ so that  the functional $F$ is bounded below and the equations above are meaningful.\\
The well-known Ljusternik-Schnirelmann principle asserts conditions on which the  $a_{(n,\alpha)}$ provide a sequence of eigenpairs  $(u_{n,\alpha},\mu_{n,\alpha})$ satisfying $\eqref{6.1}$, that is:
\begin{equation*}
			F'(u_{n,\alpha})=\mu_{n,\alpha} G'(u_{n,\alpha}).
	\end{equation*}
We express below the version corresponding to our situation.	

\begin{proposition}\label{Prop 4.1}(Ljusternik-Schnirelman principle) (see \cite{1}, \cite{4})\\
Assume that  $V > 0$ and  $(H1)-(H4)$ are fulfilled. Then  the following assertions hold.
		\begin{itemize}
								\item[(1)] If $a_{n,\alpha}>0$, then $(\ref{6.1})$ possesses a pair $\pm u_{n,\alpha}$ of eigenfunctions and an eigenvalue $\mu_{n,\alpha}\neq0$; further more $F(u_{n,\alpha})=a_{n,\alpha}$.
								\item[(2)] If $\chi_{\alpha}=\infty$, $(\ref{6.1})$ has infinitely many pairs $\pm u$ of eigenfunctions corresponding to non zero eigenvalues.
								\item[(3)] $\infty>a_{1,\alpha}\geq a_{2,\alpha}\geq...\geq0$ and $a_{n,\alpha}\rightarrow0$ as $n\rightarrow\infty$.
								\item[(4)] If $\chi_{\alpha}=\infty$ and $F(u)=0, u\in \overline{coM_{\alpha}}$ implies $\langle F'(u), u\rangle=0$, then there exists an infinite sequence $(\mu_{n,\alpha})_n$ of distinct eigenvalues of\text{ } $(\ref{6.1})$ such that $\mu_{n,\alpha}\rightarrow0$ as $n\rightarrow\infty$.
								\item[(5)] Assume that $F(u)=0, u\in \overline{coM_{\alpha}}$ implies $u=0$. Then $\chi_{\alpha}=\infty$ and there exists a sequence of eigenpairs $(u_{n,\alpha},\mu_{n,\alpha})$\text{ } of $(\ref{6.1})$ such that $u_{n,\alpha}\rightharpoonup0,\text{ } \mu_{n,\alpha}\rightarrow0$ as $n\rightarrow\infty$ and $\mu_{n,\alpha}\neq0\text{ }\forall n$.
							\end{itemize}
						\end{proposition}\text{ }\\
To apply the Ljusternik-Schnirelmann principle to our problem, we have to prove that the assumptions $(H1)-(H4)$ are satisfied.\\
Let $F$ and $G$ be defined as in relation $(\ref{2.1})$. Clearly the functionals $F$ and $G$ satisfy condition $(H1).$ As $V>0$, we obviously have $a_{n,\alpha}>0$ for any $(n,\alpha)\in\mathbb{N}^*\times\mathbb{R}^*_+$.  We next prove that $F$ and $G$ satisfy conditions $(H2)-(H4)$. For this aim, we process as follows.	
		\begin{proposition}\label{Prop4.2}(cf \cite{10,18})
	Let $\Omega$ be a domain in $\mathbb{R}^N$ and let $\phi:\Omega \times \mathbb{R}_+ \longrightarrow\mathbb{R}_+$ be a generalized  N- function which is uniformly convex and satisfies the $\Delta_2-$condition, that is there exists $c>0$ such that $\phi(x, 2t)\leq c\phi(x, t)$ for all $x \in \Omega$ and all $t\geq0$.  Then, if $(u_n)_n$ is a sequence of integrable functions in $\Omega$ such that:
							$$u(x)=\lim\limits_{n\rightarrow+\infty}u_n(x) \text{ }\text{ }\text{ for a.e }x\in\Omega \text{ }\text{ },\int_{\Omega}\phi(x, |u|)dx=\lim\limits_{n\rightarrow+\infty}\int_{\Omega}\phi\left(x, \left|u_n\right|\right)$$, we have $$\lim\limits_{n\rightarrow+\infty}\int_{\Omega}\phi\left(x, \left|u_n-u\right|\right)dx=0.$$
						\end{proposition}

						\begin{lem}\label{Lem 4.1}
							The functional $F$ satisfies condition $(H2)$.
						\end{lem}
						\begin{pre} Let $\left(u_n\right)_n\subset X$ and $u\in X$ such that
							$u_n\rightharpoonup u$. We have to show that $F'\left(u_n\right)\rightarrow F'\left(u\right)$ in $X^*$; that is $\langle F'(u_n)-F'(u),v\rangle\rightarrow0$ for any $v\in X$.\\ Let $v\in X$. Proceeding similarly as in Proposition $\eqref{Prop 2.8}$, we have
							\begin{eqnarray*}
								\left|\langle F'(u_n)-F'(u),v\rangle\right| &=& \left|\int_{\Omega}V(x)\left(|u_n|^{q(x)-2}u_n-|u|^{q(x)-2}u\right)vdx\right|\\ 
								&\leq&\int_{\Omega}V(x)\left|\left(|u_n|^{q(x)-2}u_n-|u|^{q(x)-2}u\right)v\right| dx\\
								&\leq& \left(\frac{1}{s^-}+ \frac{1}{q^{-} s'^-} + \frac{1}{q'^{-}s'^{-}}\right)\|V\|_{s(x)}\left\||u_n|^{q(x)-2}u_n-|u|^{q(x)-2}u\right\|_{s'(x)q'(x)}\|v\|_{s'(x)q(x)}\\	
&\leq& C\left(\frac{1}{s^-}+ \frac{1}{q^{-} s'^-} + \frac{1}{q'^{-}s'^{-}}\right)\|V\|_{s(x)}\left\||u_n|^{q(x)-2}u_n-|u|^{q(x)-2}u\right\|_{s'(x)q'(x)}\|v\|
							\end{eqnarray*}
where  $C$ is a positive constant due to the Sobolev embedding.  Next we have to show that $\omega_n=|u_n|^{q(x)-2}u_n\rightarrow\omega=|u|^{q(x)-2}u$ in $L^{s'(x)q'(x)}(\Omega)$. Since the  embedding $W^{1,p(x)}(\Omega)\hookrightarrow L^{s'(x)q(x)}(\Omega)$ is compact,  $u_n \rightharpoonup u$ in $X$ implies $u_n \rightarrow u$ in $L^{s'(x)q(x)}(\Omega)$. Then $u_n(x) \rightarrow u(x)$ a.e $x\in\Omega$ and $\omega_n(x) \rightarrow w(x)$ a.e in $\Omega$. Recalling again the fact that $u_n \rightarrow u$ in $L^{s'(x)q(x)}(\Omega)$, we obtain that  $\int_{\Omega}|\omega_n|^{s'(x)q'(x)}dx \longrightarrow \int_{\Omega}|\omega|^{s'(x)q'(x)}dx$. Setting:\\
							 $\phi: \Omega \times \mathbb{R}_+\longrightarrow\mathbb{R}_+$ such that $\phi(x, t) = \varphi_{s'(x)q'(x)}(t) =t^{s'(x)q'(x)}$. It is well known $\phi$  is a generalized N-function which is uniformly convex and satisfies the  the $\Delta_2$ condition. From what precedes we have $$u_n(x) \rightarrow u(x)  \mbox{ a.e }  x \in\Omega \mbox { and }
\lim\limits_{n\rightarrow\infty}\int_{\Omega}\phi(x,|\omega_n|)dx   = \int_{\Omega}\phi(x, |\omega|)dx$$.  Accordingly, by applying Proposition $\ref{Prop4.2}$, we conclude that $$ 0 =\lim\limits_{n\rightarrow\infty}\int_{\Omega}\phi(x,|\omega_n-\omega|)dx = \lim\limits_{n\rightarrow\infty}\int_{\Omega}|\omega_n-\omega|^{s'(x)q'(x)}dx=\lim\limits_{n\rightarrow\infty}\int_{\Omega}\phi_{s'(x)q'(x)}(|\omega_n-\omega|)dx$$  and then by means of the  Proposition $\ref{Prop2.6}$ we have the convergence in the sense of the norm, that is
 $$\left\||u_n|^{q(x)-2}u_n-|u|^{q(x)-2}u\right\|_{s'(x)q'(x)}= 0.$$ 
	Therefore $\langle F'(u_n)-F'(u),v\rangle\longrightarrow0$ for all $v\in X$. Hence $F'(u_n)\longrightarrow F'(u)$ in $X^*$ and thus $F$ satisfies $(H2)$. 
						\end{pre}
						
						\begin{lem}
							The functional $G'$ satisfies $(H3)$.
						\end{lem}
						\begin{pre}
							$G'$ is continuous and bounded (Proposition $\ref{Prop 2.10}$). We next show that $G'$ satisfies condition $(S_0)$.
							Let $(u_n)_n$ be a sequence in $X$ such that $u_n\rightharpoonup u$ , $G'(u_n)\rightharpoonup v_0$ and $\langle G'(u_n), u_n\rangle\longrightarrow \langle v_0,u_0\rangle$ for some $v_0\in X^*$ and $u_0\in X$. Then
							\begin{eqnarray*}
								\langle G'(u_n)-G'(u_0), u_n-u_0\rangle &=& \langle G'(u_n), u_n-u_0\rangle-\langle G'(u_0), u_n-u_0\rangle \\
								&=&\langle G'(u_n),u_n\rangle-\langle G'(u_n),u_0\rangle-\langle G'(u_0),u_n-u_0\rangle.
							\end{eqnarray*}
Accordingly, we have  $$\lim\limits_{n\rightarrow\infty}\langle G'(u_n)-G'(u_0), u_n-u_0\rangle=0$$ and 
 since $G'$ is of type $(S_+)$, we have that $u_n\longrightarrow u_0$ as $n\longrightarrow\infty$ and thus $G'$ satisfies $(H3).$ 
						\end{pre}
						
						\begin{lem}
							The functionnal $G$ satisfies $(H4)$
						\end{lem}
						\begin{pre}
							Obviously, the level set $M_\alpha$ is bounded, $\langle G'(u),u\rangle=\int_{\Omega}|\nabla u|^{p(x)}dx=\psi(u)\geq\alpha p^->0$ and $G(tu)\rightarrow\infty$ as $t\rightarrow\infty$ for any $u\in M_\alpha$. Suppose that $\inf\limits_{u\in M_\alpha}\langle G'(u),u\rangle=\inf\limits_{u\in M_\alpha}\psi(u)=0$. Thus there is $(u_n)_n\subset M_\alpha$ such that $\lim\limits_{n\rightarrow\infty}\psi(u_n)=0$. By Remark $\ref{rmq2.2}$, we have:
							$$p^-\alpha\leq\lim\limits_{n\rightarrow\infty}\psi(u_n)\leq p^+\alpha$$ that is $$p^-\alpha\leq0\leq p^+\alpha$$ which is is a contradiction. Thus $\inf\limits_{u\in M_{\alpha}}\langle G'(u),u\rangle>0$ and then $G$ satisfies $(H4)$.
						\end{pre}
						
						Thus $F$ and $G$ satisfy conditions $(H1)-(H4)$. By the Ljusternik-Schnirelman principle (Proposition $\ref{Prop 4.1}$) we conclude that: 
						\begin{proposition}\label{Thm3.2}[Existence of Ljusternik-Schnirelman sequence] (see \cite{4})\\
							For each given $\alpha>0$, there exists an eigenpair sequence $(u_{n,\alpha},\mu_{n,\alpha})_n$ obtained from the Ljusternik-Schnirelman principle such that:
							\begin{itemize}
								\item [(1)] $G(\pm u_{n,\alpha})=\alpha$ and $F(\pm u_{n,\alpha})=a_{n,\alpha}$
								\item [(2)] $\mu_{n,\alpha}=\frac{\langle F'(u_{n,\alpha}),u_{n,\alpha}\rangle}{\langle G'(u_{n,\alpha}),u_{n,\alpha}\rangle}$ where each $\mu_{n,\alpha}$ is an eigenvalue of $F'(u)=\mu G'(u)$ on $M_\alpha$
								\item [(3)] $u_{n,\alpha}\rightharpoonup0$ and $\mu_{n,\alpha}\longrightarrow0$ as $n\longrightarrow\infty$
								\item [(4)] $\infty>a_{1,\alpha}\geq a_{2,\alpha}\geq...\geq0$ and $a_{n,\alpha}\rightarrow0$ as $n\rightarrow\infty$.
							\end{itemize}	
						\end{proposition}

	\begin{rmq} \label{Rmk5.1}
							Let $\left(u_{n,\alpha},\mu_{n,\alpha}\right)$ be a sequence of eigenpairs satisfying ($\ref{6.1}$). Thus:
							\begin{equation*}
								\mu_{n,\alpha}=\frac{\langle F'(u_{n,\alpha}), u_{n,\alpha}\rangle}{\langle G'(u_{n,\alpha}), u_{n,\alpha}\rangle}=\frac{\phi(u_{n,\alpha})}{\psi(u_{n,\alpha})} \text{ }\text{ and then }\text{ }\frac{q^-a_{n,\alpha}}{p^+\alpha}\leq\mu_{n,\alpha}\leq\frac{q^+a_{n,\alpha}}{p^-\alpha}
							\end{equation*}
						\end{rmq}
			We next apply Ljusternik-Schnirelman principle in eigenvalue problem $(\ref{1.1})$ and we obtain the following results:
	\begin{proposition}\label{suite}
	\begin{itemize}
\item [(1)] Problem $(\ref{1.1})$ possesses a sequence $(\lambda_{n,\alpha})$ of eigenvalues obtained by using the (L-S) principle  and such that $\lambda_{n,\alpha}=\frac{1}{\mu_{n,\alpha}}$ for all $n$ and $\alpha>0$. Furthermore, $\lambda_{n,\alpha}\rightarrow+\infty$ as $n\rightarrow+\infty$.
\item [(2)] For any $n\in\mathbb{N}^*$ and $\alpha>0$, we have:
\begin{equation}\label{Eq 4.5}
\frac{p^-}{q^+}\frac{\alpha}{a_{n,\alpha}}\leq \lambda_{n,\alpha}\leq\frac{p^+}{q^-}\frac{\alpha}{a_{n,\alpha}}
\end{equation}
	 
	\end{itemize}
	\end{proposition}\text{ }
\begin{pre}
	\begin{itemize}
	\item [(1)]
 Let $u\in X=W^{1,p(x)}(\Omega)$ be an eigenfunction satisfying relation $(\ref{6.1})$.
		\begin{eqnarray*}
		F'(u)=\mu G'(u) &\Leftrightarrow& \int_{\Omega}V(x)|u|^{q(x)-2}uvdx=\mu\int_{\Omega}|\nabla u|^{p(x)-2}\nabla u\nabla vdx\\
		&\Leftrightarrow& \int_{\Omega}|\nabla u|^{p(x)-2}\nabla u\nabla vdx=\frac{1}{\mu}\int_{\Omega}V(x)|u|^{q(x)-2}uvdx
		\end{eqnarray*}
		By using the weak formulation $(\ref{3.6})$ we get:
			\begin{equation*}
				\lambda=\frac{1}{\mu}
			\end{equation*} Hence $\lambda_{n,\alpha}=\frac{1}{\mu_{n,\alpha}}$ for any $n\in\mathbb{N}^*$ and $\alpha>0$. From  the (L-S) principle we have $a_{n,\alpha}\rightarrow0$ as $n\rightarrow+\infty$; thus $\lambda_{n,\alpha}\rightarrow+\infty$ as $n\rightarrow+\infty$.
\item [(2)] 
Let $(\mu_{n,\alpha}, u_{n,\alpha})$ be an eigenpair satisfying $(\ref{6.1})$. It follows that $F(u_{n,\alpha})=a_{n,\alpha}$ and $G(u_{n,\alpha})=\alpha$. We then deduce from Remark $\ref{Rmk5.1}$ that
				\begin{equation*}\label{Eq4.5}
					\frac{p^-}{q^+}\frac{\alpha}{a_{n,\alpha}}\leq \lambda_{n,\alpha}\leq\frac{p^+}{q^-}\frac{\alpha}{a_{n,\alpha}}
				\end{equation*}
			
\end{itemize}
	\end{pre}
\subsubsection{\bf The smallest Lagrange multiplier and the (L-S) sequence}
\noindent

	Let $\lambda>0$ be the smallest  eigenvalue of problem $(\ref{1.1})$ obtain in Theorem $\ref{Thm4.2}$ by the Lagrange multipliers. Then there is $u_\lambda\in M_\alpha$ such that $u_\lambda$ is a (weak) solution of problem $(\ref{1.1})$ and:
						\begin{equation*}
							I_\lambda(u_\lambda)=\inf_{u\in M_\alpha}I_\lambda(u) \text{ while } F(u_\lambda)=\sup_{u\in M_\alpha}F(u)
						\end{equation*}
						
						By combining the different results in this section, we obtain the following result on the coincidence of  the  smallest eigenvalue $\lambda$ on the sphere with the first term of the (L-S) eigenvalues sequence:
						\begin{thm}\label{Thm 4.1}
							Suppose that assumption $\eqref{A}$ and $q^+ < p^-  $ hold and let $\alpha>0$. The eigenpair $(\lambda,u_\lambda)$ of problem $(\ref{1.1})$  where $\lambda $ is the smallest eigenvalue on the sphere  obtained in Theorem $\ref{Thm4.2}$ is such that:
							$$I_\lambda(u_{1,\alpha})=\inf_{u\in M_\alpha}I_\lambda(u)=I_\lambda(u_\lambda),$$ that is $u_{1,\alpha}$ is an eigenfunction solution of problem $(\ref{1.1})$  associated to the eigenvalue $\lambda$ and hence $\lambda = \lambda_{1,\alpha}= \frac{1}{\mu_{1,\alpha}}$ where $\mu_{1,\alpha}$ is the first term in the (L-S) eigenvalues sequence.
						\end{thm}

						\textbf{Proof of Theorem $\ref{Thm 4.1}:$} 
						
						Let $(\lambda,u_\lambda)$ be an eigenpair as in Theorem $\ref{Thm4.2}$. Then:
						\begin{equation}\label{Eq 4.10}
							F(u_\lambda)=\sup_{u\in M_\alpha}F(u)\geq F(u_{1,\alpha})=a_{1,\alpha}
						\end{equation}
						Put $H_0=\{\pm u_\lambda\}$. Of course, $H_0\subset M_\alpha$, compact and $\gamma(H_0)=1$ where $\gamma(H_0)$ is the genus of $H_0$ (we refer to \cite{44,38} for more details on the genus). Consequently, we have $H_0\in \Sigma_{(1,\alpha)}$ and:
						\begin{eqnarray*}
							a_{1,\alpha} &=&\sup_{H\in\Sigma_{(1,\alpha)}}\inf_{u\in H}F(u)\\
							&\geq& \inf_{u\in H_0}F(u)=F(\pm u_\lambda)=F(u_\lambda)
						\end{eqnarray*}

	This together with relation $(\ref{Eq 4.10})$ gives
		\begin{equation}\label{Eq 4.13}
		a_{1,\alpha}=\sup_{u\in M_\alpha}F(u)=F(u_{1,\alpha})=F(u_\lambda)
	\end{equation}
		Hence
		\begin{equation}
		\inf_{u\in M_\alpha}I_\lambda(u)=\alpha-\lambda\sup_{u\in M_\alpha}F(u)=\alpha-\lambda F(u_{1,\alpha})=I_\lambda(u_{1,\alpha})
		\end{equation}
		Accordingly $u_{1,\alpha}$ is a critical point of $I_\lambda$ constrained to 	$M_{\alpha}$ and then satisfied
\begin{equation}\label{comp1}
			\frac{1}{\lambda}F'(u_{1,\alpha})= G'(u_{1,\alpha}).
	\end{equation}
On the other hand $u_{(1,\alpha)}$ is an eigenfunction associated to $\lambda_{(1, \alpha)}$ and consequently we get 
\begin{equation}\label{comp2}
			\frac{1}{\lambda_{(1, \alpha)}}F'(u_{1,\alpha})= G'(u_{1,\alpha}).
	\end{equation}
From (\ref{comp1}) and (\ref{comp2}), we derive that $ \lambda = \lambda_{(1, \alpha)}$ and then the proof is complete.\\

\section{Superlinear problem}\label{section7}

Throughout this section , we will suppose  that 
\begin{equation}\label{B}
		  1<p(x)<q(x)<N < s(x) \text{ } \text { and }  s'(x)q(x) < p^{*}(x)\quad \text{ for all }x\in\overline{\Omega}
	\end{equation}
 First of all we notice that Lemma  $\ref{eqlem 3.4}$ in section $\ref{section 3}$  is expressed for any exponents $p(.)$ and $q(.)$ and accordingly  is still valid. In particular the assertions in Remark $\ref{remarq1}$ are  still holding in the following variant. 
\begin{rmq}\label{remarq2}
 Suppose $\alpha p^+ \geq 1$  then 
\begin{equation}\label{Eq 3.6'}
			\lambda_\alpha:=\alpha^{1 -\frac{q^+}{p^-}}\frac{q^- \left( p^+\right)^{-\frac{q^+}{p^-}}}{2C_HC^{q^\pm}\|V\|_{s(x)}}
		\end{equation}
and hence $\lim\limits_{\alpha \rightarrow +\infty}\lambda_\alpha = 0 \quad \text{when } \quad q^+ > p^-.$ \\
In the case that  $\alpha p^+ < 1$ 
\begin{equation}\label{Eq 3.6''}
			\lambda_\alpha:=\alpha^{1 -\frac{q^-}{p^+}}\frac{q^- \left( p^+\right)^{-\frac{q^-}{p^+}}}{2C_HC^{q^\pm}\|V\|_{s(x)}}.
		\end{equation}
Hence, $\lim\limits_{\alpha \rightarrow 0}\lambda_\alpha = +\infty$  when $q^- > p^+ $  and  if $ q^- = p^+$, then 
\begin{equation}
\forall \quad \alpha > 0, \quad \lambda_\alpha=\lambda_{p^+, q^-}: =\frac{q^- }{2p^+C_HC^{q^\pm}\|V\|_{s(x)}}.
\end{equation}
   
\end{rmq}
Here below is  the main result of the section
\begin{thm}\label{Thm 4.2}
							Suppose that assumptions $\eqref{B}$ hold  and $q^- > p^+$.
Then $\Lambda = (0, +\infty)$ .

						\end{thm}
 Under our assumption the problem  obviously  is non- coercive but   satisfies the Palais-Smale $(PS).$ 
\begin{lem}\label{Lem 3.5'}
Assume $\eqref{B}$   with $q^- \geq p^+ $. Then $ I_{\lambda}$ satisfies the Palais-Smale $(PS)$ condition.
\end{lem}
	\begin{pre}
	Let $(u_n)_{n\in\mathbb{N}}\subseteq X$  be a $(PS)$ sequence for $I_{\lambda}$, i.e., $(I_{\lambda}(u_n))_{n\in\mathbb{N}}\subset\mathbb{R}$ is bounded
	and $I'_{\lambda}(u_n)\rightarrow0$ as $n\rightarrow\infty$ , that is, there exists a positive constant $k\in\mathbb{R}$ such that
	\begin{equation}\label{MP1}
	|I_{\lambda}(u_n)|\leq k, \text{ for every } n\in\mathbb{N}.	
	\end{equation}
and 
\begin{equation}\label{MP2}
|\langle  I'_{\lambda}(u_n), v\rangle| \leq\epsilon_n\|v\| \quad \forall v \in X. 
\end{equation}
From \eqref{MP1} and  \eqref{MP2}, we have  respectively
\begin{equation}\label{MP3}
k \geq  \frac{1}{p^+}\int_{\Omega}|\nabla u_n|^{p(x)}dx-\frac{1}{q^-}\lambda\int_{\Omega}V|u_n|^{q(x)}dx \
\end{equation}
\begin{equation}\label{MP4}
 \frac{1}{q^-}\langle I'_{\lambda}(u_n), u_n\rangle  = \frac{1}{q^-}\int_{\Omega}|\nabla u_n|^{p(x)}dx-\frac{1}{q^-}\lambda\int_{\Omega}V|u_n|^{q(x)}dx\leq \frac{1}{q^-} \epsilon_n\|u_n\|
\end{equation}
Subtracting \eqref{MP2} and  \eqref{MP4}, we get 
\begin{equation}\label{MP5}
k  - \frac{1}{q^-} \epsilon_n\|u_n\| \geq  (\frac{1}{p^+} -  \frac{1}{q^-})\int_{\Omega}|\nabla u_n|^{p(x)}dx 
\end{equation}
Since $q^- \geq  p^+ $ the sequence $(u_n)_n$ is bounded in $X$ and then $I_{\lambda}$ satisfies $(PS)$ condition.\\
\end{pre}
The proof of Theorem $\ref{Thm 4.2}$ will be concluded using a standard version of the Mountain-Pass theorem (cf \cite{25}).
\begin{thm}(Mountain-Pass theorem)\\
Let $X$ a Banach space and $\Theta : X \to \mathbb{R}$ a $ C^{1}$ functional which satisfies the $(PS)$ condition. Let $S$ be a closed subset of $X$ which disconnect $X$. Let $e_0$ and $e_1$ be point of $X$ which are in distinct connected components  of $X\setminus S.$ Suppose that $\Theta$ is bounded below in $S$ and in fact the following condition is verified

\begin{equation}
\inf_{S} \Theta \geq b \textrm{ and } \max ( \Theta (e_0 ),  \Theta (e_1 )) < b.
\end{equation}
Let $$\Gamma  = \{ f \in (C([0, 1]; X),\quad  f(0) = e_0 , f(1) = e_1) \}.$$ Then 
\begin{equation}
 c = \inf_{\Gamma}\max_{t\in [0 , 1]} \Theta(f(t)) \geq b
\end{equation}
and is a critical value , that is there exists $u_0 \in X$ such that $\Theta(u_0) = c  \mbox{ and }  \Theta'(u_0) = 0)$ \\

\end{thm}
\textbf{Proof of Theorem $\ref{Thm 4.2}:$}\\
The proof consists in showing that the geometry  of the Mountain-Pass theorem is realized  with  $\Theta= I_{\lambda}.$ Since assumptions $ \eqref{B}$ and $q^- >  p^+ $ are fulfilled, then  the functional $I_{\lambda}$ satisfied the $(PS)$ condition for any $\lambda > 0 $. \\Next  for any $\alpha,$ choose $S$ to be $M_{\alpha};$ then for any $\lambda > 0, $  Lemma $\ref{eqlem 3.4}$ and the second property in Remark $\ref{remarq2}$ provide $\alpha$ (small enough) and $\lambda_{\alpha} > 0 $ such that for any $\lambda \in (0, \lambda_{\alpha}),  I_{\lambda}(u) \geq \frac{\alpha}{2}$. \\On the other hand, recalling again  assumptions $\eqref{B}$, we have $p^+<q^+$. Let $\varepsilon_0>0$ be such that $p^+ < q^+ - \varepsilon_0$. By the continuity of $q(.)$, we deduce the existence of an open set $\Omega_0\subset\Omega$ such that is $p^+ \le q^+ -\varepsilon_0\leq q(x) $ for all $x\in\Omega_0$.
		
		Let $v_0\in M_\alpha$. It is obvious that $tv_0\in X \setminus  M_{\alpha} $ for any $t > 1$. Then, we have:
		\begin{eqnarray*}
			I_\lambda(tv_0) &=& \int_{\Omega}\frac{t^{p(x)}}{p(x)}|\nabla v_0|^{p(x)}dx-\lambda\int_{\Omega}\frac{t^{q(x)}}{q(x)}V|v_0|^{q(x)}dx\\
			&\leq& t^{p^+}\alpha-\frac{\lambda}{q^+}\int_{\Omega}t^{q(x)}|v_0|^{q(x)}dx\\
			&\leq& t^{p^+}\alpha-\frac{\lambda}{q^+}\int_{\Omega_0}t^{q(x)}V(x)|v_0|^{q(x)}dx\\
			&\leq& t^{p^+}\alpha-\frac{\lambda t^{q^+-\varepsilon_0}}{q^+}\int_{\Omega_0}V(x)|v_0|^{q(x)}dx
		\end{eqnarray*}
	Therefore  for any $\lambda\in (0, \lambda_{\alpha})$ there is some $t$ large enough, say $t\ge \eta > 1$ such that 
		$I_\lambda(tv_0)<0$. Choose $e_0 = 0 $ and $e_1 = tv_0$ for a $t\ge \eta > 1$, we get $$\max(I_{\lambda}(0) = 0 , I_{\lambda}(tv_0) < 0)  < \frac{\alpha}{2} \leq \inf_{M_{\alpha}}I_{\lambda}$$ and then we are done.
\begin{cor}\label{Corol 3}
Assume that  assumptions  $\eqref{B}$ are fulfilled with $q^- = p^+ $ , then   each   $\mu \in (0, \lambda_{p^+, q^-}), $ 
 admits at least an infinitely countable and unbounded family of  eigenfunctions in $X$.
\end{cor}
\begin{pre}
Choose $\mu \in (0, \lambda_{q^+, p^-}) $ and consider an increasing  sequence of  positive real numbers  $(\alpha_n )_{n>0}$  such that $\alpha_n p^+ <  1\quad  \forall  n. $
 Under the assumptions  $\eqref{B}$ and  $q^- = p^+ $, $I_{\mu}$ satisfies $(PS)$ condition and 
\begin{equation}\label{eqCorol 3}
		I_{\mu}(u)\geq \frac{\alpha_n}{2}\quad
\textrm {  for any } u\in M_{\alpha_n}, \quad \forall n.
\end{equation}
Set $\Gamma_{\alpha_n}  = \{ f \in (C([0, 1]; X),\quad  f(0) = 0 , f(1) = t_n v_{n}^0) \}$  where $0 \in  M_{\alpha_n}$ and $t_n v_{n}^0 \in X \setminus M_{\alpha_n}$,  $v_{n}^0 \in M_{\alpha_n}$ and  $t_n $  a fixed real in  $(0, 1)$  such that  $I_{\mu}(t_n v_{n}^0) < 0$.\\ 
Moreover
 $$\max(I_{\mu}(0) = 0 , I_{\mu}(e^1_n) < 0)  < \frac{\alpha_n}{2} \leq \inf_{M_{\alpha_n}}I_{\mu}, \quad \forall n .$$
  Next, applying the Mountain-Pass Theorem, there exists at least  a sequence $(u_n)_n \in X $  of eigenfunctions and a sequence of real numbers 
$(c_n)_n$  with   $c_n = \inf\limits_{\Gamma_{\alpha_n}}\max\limits_{t\in [0 , 1]} I_{\mu}(f(t)) ,$ such that $I_{\mu}(u_n) = c_n \geq \alpha_n  \mbox{ and }  I'_{\mu}(u_n)  = 0$. Clearly  $I_{\mu}(u_n)$ tends to $+\infty$ when $n$ tends to $+\infty$ and consequently $ (u_n)_n$ is unbounded in $X.$
\end{pre}
	
{\bf Conclusion}\\
We observe that 
\begin{itemize}
\item When the exponents $p(x) \geq q(x)$, there is a neighborhood of $0$ such that any positive element is an eigenvalue and if moreover the ranges of $p(.)$ and $q(.)$ do not interfere all real number in  $(0, +\infty)$ are eigenvalues. When we constrain the eigenvalue problem to a sphere, the smallest Lagrange multiplier coincide with the first term in the (L-S) sequence.
\item When the exponents $p(x) \leq q(x)$ and when  the ranges of $p(.)$ and $q(.)$ do not interfere, we also have that all real numbers in  $(0, +\infty)$ are eigenvalues.
\item When $p(x) \geq q(x)$  and $ p^-= q^+$ or  $p(x) \leq q(x)$ and $ p^+ = q^-$, there are some eigenvalues having  an infinite sequence of eigenfunctions..
\end{itemize}
 Many challenges are still arsing from the nonhomogeneous  eigenvalues problem . For instance
\begin{itemize}
\item When $p(.) = q(.) =$ cste, the (L-S)  eigenvalues sequence  $(\lambda_{n,\alpha})_n$ given by  Proposition $\ref{suite}$ is  nondecreasing.  Is it the same in the case $p(.)$ and  $q(.)$ are functions?
\item We  still do not  know  whether the first eigenvalue in the (L-S) sequence is simple or not,  neither when $p(.) =  q(.)$  nor for $p(.) \neq  q(.)$.
\end{itemize}

\noindent
{Aboubacar Marcos $^{a}$\and  Janvier Soninhekpon $^{b}$.}\\
$^{a,b}$Institut de Math\'{e}matiques et de Sciences Physiques,
01 BP 613 Dangbo.\\
Universit\'e d'Abomey-Calavi,  B\'{e}nin.\\
 E-mail: $^{a}$ abmarcos@imsp-uac.org \\
E-mail: $^{b}$ janvier.soninhekpon@imsp-uac.org\\
\end{document}